\documentclass[preprint,12pt]{elsarticle}
\usepackage{amssymb}
\usepackage{latexsym}
\usepackage{amsmath,amsfonts,amssymb}
\usepackage{tikz}
\usepackage{hyperref}
\usepackage{graphicx}
\usepackage{color}

\newtheorem{theorem}{Theorem}[section]
\newtheorem{lemma}[theorem]{Lemma}
\newtheorem{proposition}[theorem]{Proposition}
\newtheorem{corollary}[theorem]{Corollary}
\newtheorem{remark}[theorem]{Remark}

\newtheorem{example}[theorem]{Example}
\newtheorem{definition}[theorem]{Definition}

\newcommand{\Rp}{{\mathbb R}_+}
\newcommand{\Rpm}{{\mathbb R}_+^m}
\newcommand{\Rpn}{{\mathbb R}_+^n}
\newcommand{\Rpnn}{{\mathbb R}_+^{n\times n}}
\newcommand{\Rpmn}{{\mathbb R}_+^{m\times n}}

\newcommand{\crit}{{\mathcal G}_c}
\newcommand{\digr}{{\mathcal G}}

\def\per{\mathop{\rm per}}

\def\gcd{\mathop{\rm gcd}}

\def\attr{\mathop{\rm attr}}
\def\modd{\mathop{\rm mod}}
\def\kd{\hfil$\square$\linebreak}

\newcommand{\unda}{{\underline a}}
\newcommand{\ovea}{{\overline a}}

\newcommand{\undx}{{\underline x}}
\newcommand{\ovex}{{\overline x}}

\newcommand{\Attr}{\operatorname{attr}}

\newcommand{\bfX}{\mathbf{X}}

\def\mbf#1{\mbox{\boldmath$#1$}}

\def\O{\mathcal{O}}

\def\Akr{\mathcal{D}}
\def\Dkr{\mathcal{E}}
\def\Circulant{\mathcal{Z}}
\def\IntCirculant{\mathcal{Z}^C}
\def\Nat{\mathbb{N}}
\def\per{\operatorname{per}}

\def\bfA{\mbf{A} }
\def\bfX{\mbf{X}}

\def\attr{\operatorname{attr}}

\def\N0{N\!\setminus\!\{0\}}
\def\Np{N_0}

\def\bfa{\mbf{a}}

\def\hk{\!\!\raise.5em\hbox{$^\oplus$}}
\def\kd{\hfil$\square$\linebreak}

\def\mbf#1{\mbox{\boldmath$#1$}}
\begin{document}

%
\setcounter{footnote}{1}

\begin{frontmatter}
\setcounter{footnote}{1}
 
\title{Reachability of eigenspaces for interval circulant matrices in max-algebra}
\author[rvt]{J\'an Plavka}
\ead{jan.plavka@tuke.sk}

\author[rvt2]{Serge{\u\i} Sergeev\corref{cor}\fnref{fn2}}
\ead{sergiej@gmail.com}

\address[rvt]{Department of Mathematics and Theoretical Informatics,  Technical  University,\\
B. N\v emcovej 32, 04200 Ko\v sice, Slovakia}

\address[rvt2]{University of Birmingham, School of Mathematics, Edgbaston B15 2TT, UK}

\cortext[cor]{Corresponding author.}
\fntext[fn2]{Supported by EPSRC grant EP/P019676/1}

\begin{abstract}
A nonnegative matrix $A$ is said to be strongly robust if its max-algebraic 
eigencone is universally reachable, i.e., if the orbit of 
any initial vector ends up with a max-algebraic eigenvector of $A$.  
Consider the case when the initial vector is restricted to an interval and 
$A$ can be any matrix from a given interval of 
nonnegative circulant matrices. The main aim of this paper is to classify and characterize the six types of interval robustness in this situation.
This naturally leads us also to study the max-algebraic spectral theory of circulant matrices
and the relation of inclusion between attraction cones of 
circulant matrices in max-algebra.
\end{abstract}
\end{frontmatter}



\section{Introduction}


Max-algebra has applications in such fields as discrete event systems and scheduling theory 
(among others) \cite{BCOQ, But:10, HOW}, and  plays a crucial role in the study of discrete event systems in connection with optimization problems such as scheduling or project management in which  the objective function depends on the maximum and times operations (or equivalently maximum and plus via a logarithmic transform). Notice that the main principle of discrete events  systems consisting of $n$ entities is that the entities work interactively, i.e., a given entity must wait before proceeding to its next job until
certain others have completed their current jobs. The steady states of such systems correspond to the 
max-algebraic eigenvectors of the matrices that describe them, therefore the investigation of reachability of the set of eigenvectors from a given state by a given system is important for such applications. Matrices for which the steady states of the corresponding systems
are reached with any nontrivial starting vector are called robust, see \cite{But:10} Section~8.6.

In practice, matrix entry values are not exact numbers and usually are contained within intervals, and 
therefore interval arithmetic is an efficient way to represent matrices in a guaranteed way on a computer.
A max-algebraic (tropical) version of interval analysis was developed, e.g., in
\cite{LitSob-2001}, which emphasized the polynomiality of some algorithms of max-algebraic interval analysis. That polynomiality was in 
striking contrast with NP-hardness of relevant algorithms previously known in usual interval analysis.
Independently,~\cite{fnrrz} developed a theory of some max-algebraic linear systems with interval coefficients and optimization problems over such systems. 

When developing interval extensions of linear algebra problems a whole range of solvability problems routinely arises, by considering all possible combinations of quantifiers (as in Definition~\ref{def:introb} of the present paper). 
In classical linear algebra this leads to the notions of united solutions, controllable solutions and tolerable solutions~\cite{ro2003, Shary2002}. In max-algebra we similarly have, e.g.,
four types of interval extensions of the max-algebraic
spectral problem~\cite{gaplpo01} or two types of interval extensions of robustness studied in~\cite{mp}. 

Similarly to~\cite{mp}, the present paper also considers max-algebraic interval extensions of robustness and reachability problems.  However, we 
focus on matrices of a certain special type: circulants.  In usual algebra, circulant matrices have a number of geometric applications~\cite{Davis}.
A more recent application of circulants can be found in~\cite{circ}. There, an algebraic construction based on circulant matrices allows for designing LDPC codes with efficient 
encoder implementation, in contrast to designing LDPC codes based on random construction techniques
which make it difficult to store and assess a large parity-check matrix or to analyze the performance of the
code.
 In max-algebra, circulant matrices appear to describe the periodic regime of sequences of matrix powers~\cite{But:10,Ser-11}. It is also easy to see that 
circulant matrices of a given dimension form a commutative semigroup, both in max-algebra and in usual linear algebra.

When considering matrices of special
type, it is natural to require that the set of matrices that is an interval extension of such a matrix 
can contain matrices of that type only. This is a basic idea behind the notion of interval circulant matrix
defined here. The main aim of the present paper is thus to 
classify and characterize the six types of interval robustness for circulant matrices in max-algebra. However, 
obtaining such a characterization is not possible without a deeper study of properties of 
circulant matrices in max-algebra, which is itself of some theoretical interest.

We now outline the organization of the paper and the
results obtained there. Section~\ref{preli} is devoted to some basic notions of max-algebra and its connections to the theory of digraphs and max-algebraic  
convexity. In particular, we revisit the max-algebraic spectral theory here, focusing on the
eigencone and the attraction cone associated with an arbitrary eigenvalue, 
the cyclicity of critical graphs and the ultimate
periodicity of max-algebraic matrix powers and orbits.

Section~\ref{s:circulants} presents some known as well as some new results on the 
spectral theory and attraction cones of circulant matrices. In particular,
Proposition~\ref{t:critcirc} describes the critical node sets of circulant matrices and presents 
several formulae for the cyclicity of the critical graph of a circulant matrix. 
This result combines together some facts that have been previously obtained 
or stated in~\cite{mmp,mp,Tom-10}. The main new result of this section is Theorem~\ref{t:attrincl}, which deals with a particular problem of inclusion of the attraction cones of circulant matrices $A$ and $B$ satisfying $A\leq B$ and having the same maximum cycle mean. It appears that inclusion $\attr(A)\subseteq\attr(B)$ holds for such 
circulant matrices. Note that it does not hold for general matrices, as Example~\ref{ex:counter} demonstrates.
Section~\ref{s:circulants} also contains several motivating examples. The proofs of Proposition~\ref{t:critcirc} and Theorem~\ref{t:attrincl} are deferred to Section~\ref{s:proofs}.

Based on the result about inclusion of attraction cones of Theorem~\ref{t:attrincl},
Section~\ref{s:interval} characterizes various types of
interval robustness which are described in Definition~\ref{def:introb}. 
Some of them can be verified in polynomial time, see Theorems~\ref{t:posrob},~\ref{t:unirob},~\ref{T_toleranceA_robust}.  Other types of robustness reduce to 
max-algebraic two-sided systems of equations and inequalities for which efficient algorithms
exist but the problem of constructing a polynomial algorithm remains open. See Theorems~\ref{TXRCIR},
~\ref{T_weak_tolerance_robust},~\ref{T_possiblyA_robust}. 


Subsection~\ref{ss:proof1} presents a proof of Proposition~\ref{t:critcirc}.
The proof uses the fact that any circulant matrix is strictly visualized in the 
sense of~\cite{SSB} and relies in part on the results of~\cite{Gav-99,Gav:04}.

Subsection~\ref{ss:proof2} presents a proof of Theorem~\ref{t:attrincl}. In particular,
the proof draws upon the role of cyclic classes in the max-linear systems of equations describing
attraction cones, as presented in~\cite{But:10}~Chapter~8 and~\cite{Ser-11}.

\section{Preliminaries}\label{preli}

\subsection{Main definitions and problem statements}

By max-algebra we mean the set of nonnegative numbers $\Rp$ equipped with the 
usual multiplication $a\cdot b$ and the idempotent addition $a\oplus b:=\max(a,b)$. 
These arithmetical operations are then routinely extended to matrices and vectors:
in particular, $(A\otimes B)_{i,k}=\bigoplus_{j} A_{i,j}\cdot B_{j,k}$ and 
$(A\oplus B)_{i,j}=A_{i,j}\oplus B_{i,j}$ for any two nonnegative matrices $A,B$
of appropriate sizes. We will also consider the max-algebraic powers of matrices
$A^k:=\underbrace{A\otimes\ldots\otimes A}_{k}.$

In what follows, we will be interested in the orbits of vectors 
under the action of matrices, that is, the sets
\begin{equation}
\label{e:orbit}
\O(A,x)=\{x,\, A\otimes x,\, {A^2\otimes x},\ldots\},
\end{equation}
and especially in the case when the orbit of a vector hits an eigenvector of $A$.
Let us now give formal definitions related to the max-algebraic eigenproblem.

\begin{definition}[Eigenvalues and Eigenvectors]
\label{def:eigenproblem}

{\rm A value $\lambda\in\Rp$ is called a {\em (max-algebraic) eigenvalue}
of $A\in\Rpnn$ if $A\otimes x=\lambda x$ for
some $x\in\Rpn\backslash\{0\}$. The greatest eigenvalue of $A$ will be denoted by $\lambda(A)$.

A vector $x\in\Rpn\backslash\{0\}$ satisfying
$A\otimes x=\lambda x$ is called a {\em (max-algebraic) eigenvector}
associated with $A$.

The {\em eigencone of $A$} associated with 
eigenvalue $\lambda$ is defined as the set containing 
all eigenvectors of $A$ with associated eigenvalue $\lambda$ as well as the zero vector:
$$V(A,\lambda)=\{x\in \Rpn \colon \ A\otimes x= \lambda\otimes x\}.$$
}
\end{definition}

One of the key notions of the paper is that of attraction cone:
the set which comprises all vectors whose orbit hits a given eigencone.

\begin{definition}[Attraction cones]
\label{def:attrcone}
{\rm The {\em attraction cone} of $A\in\Rpnn$ associated with eigenvalue $\lambda$ is the set
$$\attr(A,\lambda) = \{x \in \Rpn\colon \ \
\O(A, x) \cap V(A,\lambda)\neq \emptyset\}.$$

We also denote $\attr(A)=\attr(A,\lambda(A))$.}
\end{definition}

Any eigencone or any attraction cone is a max cone, in the sense of the 
following definition.

\begin{definition}[Max cones]
\label{def:maxcone}
{\rm 
A set $V\subseteq\Rpn$ is called a 
{\em max cone} if for all $x\in V$, $y\in V$ any {\em max-linear combination} 
$\alpha x\oplus\beta y$ (where $\alpha,\beta\in\Rp$) belongs to $V$.}
\end{definition}

We will use the following notational shortcuts. 

\begin{definition}[Index Sets $N$ and $N_0$]
\label{def:NN0}
{\rm We denote}
$$N=\{1,\ldots,n\},\quad N_0=\{0,\ldots,n-1\}.$$ 
\end{definition}

In this paper we deal with the following special class of matrices in max-algebra.
\begin{definition}[Circulant Matrices]
\label{def:circulant} 
{\rm A matrix $A\in \Rpnn$ is called circulant, if it has entries $A_{i,j}=a_t$ 
for $i,j\in N$, $t\in N_0$ such that $t\equiv (j-i)(\modd\; n)$ and  
$a_0,a_1,\ldots,a_{n-1}\in\Rp$. Equivalently, $A$ is a circulant matrix if it is of the form
	\begin{displaymath}
	A=\left(\begin{array}{cccccc}
	a_0&a_1&a_2&\dots&a_{n-2}&a_{n-1}\\
	a_{n-1}&a_0&a_1&\dots&a_{n-3}&a_{n-2}\\
	\vdots&\vdots&\vdots&\ &\vdots&\vdots\\
	a_1&a_2&a_3&\dots&a_{n-1}&a_0\\
	\end{array}\right).
	\end{displaymath} 
	for some $a_0,a_1,\ldots,a_{n-1}\in\Rp$.
Such a circulant matrix 
we will also denote by $\Circulant(a_0,\dots,a_{n-1})$.}
\end{definition}

Circulant matrices will be the main topic of Section~\ref{s:circulants} and Section~\ref{s:proofs}, where 
we will study their spectral theory and attraction cones.

\if{
The robustness of a matrix $A$ means that the members of $\O(A,x)$ achieve  the eigenspace of $A$ starting at arbitrary initial vector $x\in \Rpn$. In the following we shall deal with the so-called $\bfX$-robustness of matrix which requires achieving the eigenspace  starting at each  vector $x$ from a given  interval vector $\bfX$ provided that $\O(A,x)\subseteq \bfX$.
}\fi

The final part of this paper is devoted to intervals and interval circulant matrices.

\begin{definition}[Intervals]
\label{def:intervals}
{\rm A set $\mbf{X}\subseteq\Rpn$ is called an {\em interval} if it is of the form
$$\mbf{X}=\times_{i=1}^n \mbf{X}_i,$$
for $\mbf{X}_i$ nonempty subsets of $\Rp$ 
taking any of the following four forms:
$$[\undx_i,\ovex_i],\ (\undx_i,\ovex_i),\ (\undx_i,\ovex_i],\ [\undx_i,\ovex_i),$$
for $\undx_i,\ovex_i\in\Rp$.}
\end{definition}


\begin{definition}[Interval Circulant Matrices]  
\label{def:intcirculant}
{\rm By $\IntCirculant(\bfa_0,\dots,\bfa_{n-1})$ we denote the 
set of all circulant matrices $A$ such that $A_{i,j}\in \bfa_t$
for $i,j\in N$ and $t\in N_0$ such that $t\equiv (j-i)(\modd\; n)$, 
where $\bfa_0,\dots,\bfa_{n-1}$ are intervals 
independently taking any of the four forms listed in Definition~\ref{def:intervals}.

A set of circulant matrices that is of the form $\IntCirculant(\bfa_0,\dots,\bfa_{n-1})$ 
for intervals $\bfa_0,\dots,\bfa_{n-1}$ is called an {\rm interval circulant matrix}}.
\end{definition}

In the literature on max-algebra, 
$A\in\Rpnn$ is called robust if $\Attr(A)=\Rpn$, see~\cite{But:10} Section 8.6.
In this paper we consider various extensions of this notion 
to interval circulant matrices.
These extensions are listed in the following definition.

\begin{definition}[Interval Robustness]
\label{def:introb}
	{\rm Let $\mbf{X}\subseteq\Rpn$ be an interval and 
	$\IntCirculant(\bfa_0,\dots,\bfa_{n-1})$ be an interval circulant matrix. Then $\IntCirculant(\bfa_0,\dots,\bfa_{n-1})$  is called
	\begin{itemize}
		\item\hspace{-1.1cm}
		{\em possibly $\bfX-$robust}  if
		$(\exists A\in \IntCirculant(\bfa_0,\dots,\bfa_{n-1})) (\forall x\in \bfX) [\, x\in\Attr(A) \,] \enspace,$
		\item\hspace{-1.1cm}
		{\em universally $\bfX-$robust}   if
		$(\forall A\in \IntCirculant(\bfa_0,\dots,\bfa_{n-1})) (\forall x\in \bfX) [\, x\in\Attr(A) \,] \enspace,$
		\item\hspace{-1.1cm}
		{\em tolerance $\bfX-$robust}   if
		$(\forall A\in \IntCirculant(\bfa_0,\dots,\bfa_{n-1})) (\exists x\in \bfX) [\,x\in\Attr(A) \,] \enspace,$
		\item\hspace{-1.1cm}
		{\em weakly tolerance $\bfX-$robust} if
		$(\exists A\in \IntCirculant(\bfa_0,\dots,\bfa_{n-1})) (\exists x\in \bfX) [x\in\Attr(A)]$
	\end{itemize}
	and
	$\bfX$  is called
	\begin{itemize}
		\item\hspace{-1.1cm}
		{\em possibly $\IntCirculant-$robust}   if
		$(\exists x\in \bfX) (\forall A\in \IntCirculant(\bfa_0,\dots,\bfa_{n-1})) [\, x\in\Attr(A)  \,] \enspace,$
		\item\hspace{-1.1cm}
		{\em tolerance $\IntCirculant-$robust}   if
		$(\forall x\in \bfX) (\exists A\in \IntCirculant(\bfa_0,\dots,\bfa_{n-1})) [\, x\in\Attr(A) \,] \enspace.$
		
	\end{itemize}
	}
\end{definition}

In particular, the aim of Section~\ref{s:interval} will be to derive an efficient
characterization of these types of
interval robustness.

\subsection{Associated graphs, critical graphs and periodicity}

Let us start with the following basic definition. For relevant definitions see also, e.g.,
\cite{But:10} Section 1.5. 

\begin{definition}[Digraphs, Walks, Cycles and Connectivity]
\label{def:graphs}
{\rm Let $\digr$ be a digraph with set of nodes $N$ and set of edges $E$.
A {\em walk} on $\digr$ is a sequence $W=(i_0,i_1,\ldots, i_l)$ with 
$i_0,i_1,\ldots,i_l\in N$ where each pair $(i_{s-1},i_s)$ for
$s\in\{1,\ldots,l\}$ is an edge. If $i_0=i$ and $i_l=j$ then $W$ is said to be 
{\rm connecting $i$ to $j$}, and $l$ is called the {\rm length} of $W$.

$\digr$ is called {\em strongly connected} if for each 
$i,j\in N$ with $i\neq j$ there exists a walk on $\digr$ connecting $i$ to $j$.

For $A\in\Rpnn$, the {\em weighted digraph $\digr(A)$ associated with $A$}
is the digraph with set of nodes $N=\{1,\ldots,n\}$ 
and set of edges $E=\{(i,j)\colon A_{i,j}\neq 0\}$, where 
$A_{i,j}$ is the {\em weight} of an edge $(i,j)$.

If $\digr=\digr(A)$ then the {\em weight} of $W=(i_0,i_1,\ldots, i_l)$ is 
defined by $A_{i_0,i_1}\cdot A_{i_1,i_2}\cdot\ldots\cdot A_{i_{l-1},i_l}$.
This walk is called a {\em cycle} if 
$i_l=i_0$, with the {\em cycle (geometric) mean} defined by $(A_{i_0,i_1}\cdot A_{i_1,i_2}\cdot\ldots\cdot A_{i_{l-1},i_0})^{1/l}$.}
\end{definition}

Let us also give a separate definition of the maximum cycle mean. 

\begin{definition}[Maximum cycle (geometric) mean]
\label{def:mcm}
{\rm The {\em maximum cycle (geometric) mean} of any $A\in\Rpnn$ or of $\digr(A)$ is}
\begin{equation}
\label{e:mcm}
\max\limits_{k=1}^n\max\limits_{1\leq i_1,\ldots,i_k\leq n}
(A_{i_1,i_2}\cdot A_{i_2,i_3}\dots A_{i_k,i_1})^{1/k}.
\end{equation}
\end{definition}
The striking importance of this concept in max-algebra is due to the 
following fact.
\begin{proposition}[e.g., \cite{But:10}, Corollary~4.5.6]
\label{p:mcm} 
For any $A\in\Rpnn$, its greatest max-algebraic eigenvalue 
($\lambda(A)$) is equal to~\eqref{e:mcm}.
\end{proposition}


The concept of irreducible matrix is common for 
max-algebra and nonnegative linear algebra, and it is most conveniently 
defined via the associated digraph.

\begin{definition}[Irreducible, Reducible and Completely Reducible]
{\rm $A$ is called {\em irreducible} if $\digr(A)$ is strongly connected, and {\em reducible} otherwise. 

Digraph $\digr$ is called {\em completely 
reducible} if it consists of several strongly connected subgraphs called {\em components} such that there
are no walks connecting a node from one component to a node of another component.
$A$ is called {\em completely 
reducible} if so is $\digr(A)$.}
\end{definition}

Note that any irreducible matrix is completely reducible.
Observe also the following criterion of complete reducibility.

\begin{proposition}
\label{p:compred}
A digraph $\digr=(N,E)$ is completely reducible 
if and only if every edge of $E$
lies in a cycle of $\digr$.
\end{proposition}
{\it Proof.}
{\bf ``If":} Suppose that $\digr$ contains two maximal strongly connected
subgraphs $\digr_1$ and $\digr_2$ and that there is a walk connecting one subgraph to the other. Without loss of generality
we can assume that the walk does not contain nodes from any other subgraphs, so that it contains an 
edge $(i,j)$ with $i\in\digr_1$ and $j\in\digr_2$. As this edge is on a cycle, there is also a walk from
$j$ to $i$. However, this implies that $\digr_1$ and $\digr_2$ both belong to a larger strongly connected 
subgraph of $\digr$ thus contradicting their maximality. Thus the ``if" part is proved. 

{\bf ``Only if":} If $\digr$ is completely reducible then each edge 
$(i,j)$ belongs to a strongly connected subgraph of $\digr$, and it belongs to 
a cycle since there exists a walk connecting $j$ back to $i$.
\kd

The following subdigraph of $\digr(A)$ is crucial for the max-algebraic spectral theory
and it is an example of completely reducible digraph.

\begin{definition}[Critical Digraphs]
\label{def:critical}
{\rm The {\em critical digraph} of $A$, denoted by $\crit(A)$, consists of 
all nodes and edges of the cycles of $\digr(A)$ at which the maximum 
cycle mean of $A$~\eqref{e:mcm} is attained. These cycles are called {\em critical cycles}. 
The nodes of $\crit(A)$ are called 
{\em critical nodes} and their set is denoted by $N_c(A)$, and the edges of $\crit(A)$ are called 
{\em critical edges} and their set is denoted by $E_c(A)$.}
\end{definition}

\begin{corollary}
\label{c:compred-crit}
Any critical graph is completely reducible.
\end{corollary}
{\it Proof.} By Definition~\ref{def:critical}, every edge of $\crit(A)$ belongs to a cycle of $\crit(A)$.
The claim now follows from Proposition~\ref{p:compred}.
\kd

The concept of the digraph's cyclicity is crucial for the study of attraction cones (Definition~\ref{def:attrcone}) and the ultimate periodicity of
$\{A^t\}_{t\geq 1}$ (to be defined soon).

\begin{definition}[Cyclicity]
\label{def:cyclicity}
{\rm For a strongly connected digraph, its {\em cyclicity} is defined as 
the g.c.d. of the lengths of all cycles of that digraph.

Cyclicity of a completely reducible digraph 
is defined as the l.c.m. of the cyclicities of its components.

Cyclicity of a digraph $\digr$ is denoted by $\sigma(\digr)$.} 
\end{definition}

We now discuss the ultimate periodicity of max-algebraic matrix powers.

\begin{definition}[Ultimate Periodicity]
\label{def:ultperiod}
{\rm Let $\{\alpha_k\}_{k\geq 1}$ be a sequence of some elements. If 
there exists $T$ such that $\alpha_{t+\sigma}=\alpha_t$ for all $t\geq T$ and some $\sigma$ 
(i.e., $\alpha_{t+\sigma}$ and $\alpha_t$ are identical), then $\{\alpha_k\}_{k\geq 1}$ is called
{\em ultimately periodic}. The least $T$ and the least $\sigma$ for which the above property holds are called
the {\em transient} and the {\em ultimate period} of $\{\alpha_k\}_{k\geq 1}$ respectively.}
\end{definition}

\begin{proposition}[\cite{Cohen+}]
\label{p:period}
Let $A\in\Rpnn$ be an irreducible matrix with 
$\lambda(A)\neq 0$.
Then $\{(A/\lambda(A))^t\}_{t\geq 1}$ is ultimately periodic and $\sigma(\crit(A))$ is the ultimate period
of that sequence.
\end{proposition}

In this paper we also need the following trivial extension of Proposition~\ref{p:period}
and its consequence for orbits of vectors.

\begin{corollary}
\label{c:compred}
Let $A\in\Rpnn$ be a completely reducible matrix with $\lambda(A)\neq 0$, such that the maximum
cycle mean of each component of $\digr(A)$ is the same (and equal to 
$\lambda(A)$).
Then $\{(A/\lambda(A))^t\}_{t\geq 1}$ is ultimately periodic and $\sigma(\crit(A))$ is the ultimate period
of that sequence.
\end{corollary}

\begin{corollary}
\label{c:period}
Under the conditions of Proposition~\ref{p:period} or Corollary~\ref{c:compred},
$\{(A/\lambda(A))^t\otimes x\}_{t\geq 1}$ is ultimately periodic 
for any $x\in\Rpn$.
\end{corollary}

Let us now introduce some notation related to the ultimate periodicity.

\begin{definition}
\label{def:period}
{\rm Let $A\in\Rpnn$ have $\lambda(A)\neq 0$.
If $\{(A/\lambda(A))^t\}_{t\geq 1}$ is ultimately periodic then denote by
$T(A)$ the transient and by $\per(A)$ the ultimate period of that sequence.}
\end{definition}

Thus $\per(A)=\sigma(\crit(A))$ for any $A$ satisfying the 
condition of Proposition~\ref{p:period}
or Corollary~\ref{c:compred}.

The ultimate period of $\{(A/\lambda(A))^t\otimes x\}_{t\geq 1}$
does not necessarily equal the cyclicity of $\crit(A)$, and the 
attraction cone associated with $\lambda(A)$ consists of the vectors for which the ultimate period 
of $\{(A/\lambda(A))^t\otimes x\}_{t\geq 1}$ is equal to $1$. More precisely, 
we have the following.

\begin{proposition}
\label{p:attrdesc}
Let $A\in\Rpnn$ be a completely reducible matrix with $\lambda(A)\neq 0$ 
such that the maximum cycle mean of each component of $\digr(A)$ is the same (and equal to 
$\lambda(A)$). Then
\begin{equation*}
\Attr(A)=\{x\in\Rpn\colon \lambda(A)A^t\otimes x=A^{t+1}\otimes x\},\quad\text{where $t\geq T(A)$.}
\end{equation*}
\end{proposition}
{\it Proof.}  
By definition $x\in\Attr(A)$ if and only if 
$A^{s+1}\otimes x= \lambda(A) A^s\otimes x$ for some $s$, hence $\lambda(A)A^t\otimes x=A^{t+1}\otimes x$ for 
some $t\geq T(A)$ is sufficient for $x\in\attr(A)$. 
\if{
and
$A^s\otimes x\neq 0$ for some $s$. 
For a completely reducible $A$ we have $A^t\otimes x\neq 0$ for any 
$x\neq 0$, hence $\lambda(A) A^t\otimes x=A^{t+1}\otimes x$  is sufficient for $x\in\Attr(A)$. 
}\fi
For the necessity observe that  $A^{s+1}\otimes x= \lambda(A) A^s\otimes x$ implies $A^{s'+1}\otimes x= \lambda(A) A^{s'}\otimes x$ for some $s'\geq \max(s,T(A))$ and such that 
$(A/\lambda(A))^{s'}=(A/\lambda(A))^{t}$, and hence $A^{t+1}\otimes x= \lambda(A) A^t\otimes x$.
\kd

\begin{corollary}
\label{c:attrclosed}
 Under the conditions of Proposition~\ref{p:attrdesc} $\attr(A)$ is a closed max-cone.
\end{corollary}
 {\it Proof.}
Under these conditions $\Attr(A)$ is the solution set of the system 
$\lambda(A)A^t\otimes x=A^{t+1}\otimes x$. This solution set is a max-cone since it is closed under
taking max-linear combinations (see Definition~\ref{def:maxcone}) and it is a closed set since 
all arithmetic operations of max-algebra are continuous.
\kd

Let us finally consider the attraction cones of the following two matrices satisfying the 
conditions of Proposition~\ref{p:attrdesc}.

\begin{example}
\label{ex:counter}
{\rm Take
\begin{equation*}
A=
\begin{pmatrix}
0.5 & 1 & 0.2 & 0\\
1   & 0.5 & 0.2 & 0\\
0.2 & 0.2 & 0.2 & 0\\
0  & 0 & 0 & 1
\end{pmatrix},
B=
\begin{pmatrix}
0.5 & 1 & 0.2 & 0\\
1 & 0.5 & 0.3 & 0\\
0.4 & 0.4 & 0.4 & 0\\
0 & 0 & 0 & 1
\end{pmatrix}
\end{equation*}
The ultimate periods of $\{A,A^2,A^3,\ldots\}$
and $\{B,B^2,B^3,\ldots\}$ equal $2$. In the first case, the periodicity 
starts from $A^2$ (i.e., we have $A^2=A^4$), and in the second case it starts from 
$B^3$ (i.e., we have $B^3=B^5$). The attraction cones are
\begin{equation*}
\Attr(A)=\{x\colon A^3\otimes x=A^4\otimes x\},\quad 
\Attr(B)=\{x\colon B^3\otimes x=B^4\otimes x\},
\end{equation*}
where
\begin{equation*}
\begin{split}
& A^3=
\begin{pmatrix}
0.5 & 1 & 0.2 & 0\\
1 & 0.5 & 0.2 & 0\\
0.2 & 0.2 & 0.04 & 0\\
0 & 0 & 0 & 1 
\end{pmatrix},\quad
A^4=
\begin{pmatrix}
1 & 0.5 & 0.2 & 0\\
0.5 & 1 & 0.2 & 0\\
0.2 & 0.2 & 0.04 & 0\\
0 & 0 & 0 & 1
\end{pmatrix},\\
& B^3=
\begin{pmatrix}
0.5 & 1 & 0.2& 0\\
1 & 0.5 & 0.3 & 0\\
0.4 & 0.4 & 0.12 & 0\\
0 & 0 & 0 & 1
\end{pmatrix},\quad
B^4=
\begin{pmatrix}
1 & 0.5 & 0.3 & 0\\
0.5 & 1 & 0.2 & 0\\
0.4 & 0.4 & 0.12 & 0\\
0 & 0 & 0 & 1
\end{pmatrix}.
\end{split}
\end{equation*}
We further see that in both cases, the systems defining these attraction cones reduce to just 
one equation:
\begin{equation*}
\begin{split}
\Attr(A)&=\{x\colon 0.5 x_1\oplus x_2\oplus 0.2 x_3= x_1\oplus 0.5x_2\oplus 0.2 x_3\},\\ 
\Attr(B)&=\{x\colon 0.5 x_1\oplus x_2\oplus 0.2 x_3= x_1\oplus 0.5x_2\oplus 0.3 x_3\},
\end{split}
\end{equation*}
Observe that $x=[1\;1\;5\; 1]$ belongs to $\Attr(A)$ but not to $\Attr(B)$, 
and  $x=[0.5\;1\;\frac{10}{3}\; 1]$ belongs to $\attr(B)$ but not to $\attr(A)$.\\
 }\kd
\end{example}

Example~\ref{ex:counter} also shows that Theorem~\ref{t:attrincl}, the main result of the next section 
which claims that $\attr(A)\subseteq\attr(B)$ for any circulant $A,B$ with $A\leq B$ and 
$\lambda(A)=\lambda(B)$, is not true for general completely reducible (or irreducible) 
matrices.

\section{Circulant matrices: critical graph and attraction cones}
\label{s:circulants}

Let us start with the following statement, which is well known in usual linear algebra. See, e.g.,
\cite{Davis} Theorem 3.1.1. A proof of it in max-algebra, 
which works equally well in the usual linear algebra case, is given below for 
the reader's convenience. 
  
\begin{proposition}
\label{p:classical}
Let $A,B\in\Rpnn$ be circulant matrices. Then $A\otimes B$ is also 
circulant. In particular, any max-algebraic power of $A$ (or $B$) is a circulant.
\end{proposition}
{\it Proof.} Observe that $A$ is a circulant matrix if and only if we can represent 
$A=a_0I\oplus a_1P\oplus\ldots\oplus a_{n-1}P^{n-1},$ where
\begin{equation*}
P=
\begin{pmatrix}
0 & 1 & 0 &\ldots & 0\\
0 & 0 & 1 &\ldots & 0\\
\vdots &\ldots &\ddots &\ldots &\vdots\\
0 & 0 & \ldots &  0&  1\\
1 & 0 & \ldots & 0 & 0
\end{pmatrix}.
\end{equation*}
In this case $A=\Circulant(a_0,a_1,\ldots,a_{n-1})$. 
Computing $A\otimes B$ amounts to multiplying $a_0I\oplus a_1P\oplus\ldots\oplus a_{n-1}P^{n-1}$
by $b_0I\oplus b_1P\oplus\ldots\oplus b_{n-1}P^{n-1}$, assuming that 
$A=\Circulant(a_0,a_1,\ldots,a_{n-1})$ and $B=\Circulant(b_0,b_1,\ldots,b_{n-1})$.
This multiplication results in an expression of the form $c_0I\oplus c_1P\oplus\ldots\oplus c_{n-1}P^{n-1},$
thus also a circulant. 

Writing $A^t$ as $A^{t-1}\otimes A$ for every $t\geq 2$, we also show that $A^t$ is a circulant by 
a simple inductive argument. 
\kd

The following observation will play a key role in proving many properties of circulants.

\begin{lemma}
\label{l:elementary}
Let $A\in \Rpnn$ be a nonzero circulant matrix and let $A_{i,j}=\mu\neq 0$ for some $i,j\in N$.
Then $(i,j)$ belongs to a cycle $(i_1,\ldots,i_n,i_1)$ with 
$i_{t}-i_{t-1}\equiv (j-i)(\modd\; n)$ for all $t\in\{2,\ldots,n\}$, and
$i_1-i_n\equiv (j-i)(\modd\; n)$. The weight of each edge in $(i_1,\ldots,i_n,i_1)$ equals $A_{i,j}=\mu$.
\end{lemma}
 {\it Proof.} 
Consider an infinite sequence $\{i_{\ell}\}_{\ell\geq 1}$ where $i_{1}=i$, $i_2=j$, 
$i_{\ell+1}-i_{\ell}\equiv (j-i)(\modd\; n)$ for all $\ell\geq 1$ and $i_{\ell}\in N$ for all $\ell\geq 1$. 
By Definition~\ref{def:circulant},
$A_{i_{\ell},i_{\ell+1}}=A_{i,j}$ for all $\ell\geq 1$. However, we also have that $i_{n+1}=i_1$
since $i_{n+1}-i_1\equiv n\cdot (j-i) (\modd\; n)=0(\modd\; n)$. Hence the claim follows.
\kd


\begin{proposition}
\label{p:lambdacirc}
Let $A=\Circulant(a_0,\dots,a_{n-1})$. 
Then $A$ has a unique max-algebraic eigenvalue equal to
\begin{equation}
\label{e:lambdaAcirc}
\lambda(A)=\max_{k=0}^{n-1} a_k.
\end{equation}
If $A\neq 0$ then $\lambda(A)\neq 0$ and
all nodes in $N$ are critical.
\end{proposition}
{\it Proof.}
If $A\neq 0$ then $\max(a_0,\ldots,a_{n-1})>0$. In this case, let $i$ and $j$ be such that
$A_{i,j}=\mu>0$. By Lemma~\ref{l:elementary} $(i,j)$ belongs to a cycle $(i_1,\ldots,i_n,i_1)$
where the weights of all edges are equal to $\mu$. It follows that the cycle mean of that cycle is also $\mu$.
Thus, the maximal cycle mean is equal to the 
maximal weight of edges, which shows~\eqref{e:lambdaAcirc}. Taking $k$
such that $a_k=\lambda(A)$, for each $i\in N$ we have $j$ with $k\equiv (j-i)(\modd\; n)$ such that
$A_{i,j}=\lambda(A)$, hence each $i\in N$ is on a critical cycle.  Since all nodes $\digr(A)$ are critical, 
$A$ has a unique eigenvalue equal to $\lambda(A)$ as it follows, e.g., from \cite{But:10} Corollary 4.5.8.

If $A=0$ then $\max(a_0,\ldots,a_{n-1})=0=\lambda(A)$.
\kd

Note that
equation~\eqref{e:lambdaAcirc} was obtained already in~\cite{Pla-01}, Theorem 2.1. However, we preferred 
to give a partially self-contained proof of this equation for the reader's convenience. 
\begin{corollary}
	\label{l:zerocase}
	Let $A\in\Rpnn$ be a circulant matrix. Then $\lambda(A)=0$ if and only if $A=0$.
\end{corollary}
{\it Proof.} Obviously, $\lambda(A)=0$ if $A=0$. The "only if" part is equivalent 
to the implication $(A\neq 0)\Rightarrow (\lambda(A)\neq 0)$ stated in Proposition~\ref{p:lambdacirc}.
\kd

\if{
By~\eqref{e:lambdaAcirc}, $\lambda(A)$ equals the greatest entry of $A$. Therefore if $\lambda(A)=0$ then
each entry of $A$ is $0$, thus $A=0$. The remaining ``if" part is obvious.
\kd
}\fi

We now formulate the following immediate corollary of Proposition~\ref{p:compred}.

\begin{corollary}
\label{c:compred-circ}
Any circulant matrix $A$ is completely reducible.
\end{corollary}
{\it Proof.} If $A=0$ then $\digr(A)$ has no edges and is completely reducible. 
Otherwise, by Lemma~\ref{l:elementary} any edge of $\digr(A)$ belongs to a cycle, and  
 the claim follows from Proposition~\ref{p:compred}. 
\kd

\begin{proposition}
\label{p:transient-circ}
For any nonzero circulant matrix $A\in\Rpnn$ the matrix sequence $\{(A/\lambda(A))^t\}_{t\geq 1}$
is ultimately periodic, and $T(A)\leq (n-1)^2+1$.
\end{proposition}
{\it Proof.}
For the first part of the claim observe that any circulant matrix is
completely reducible by Corollary~\ref{c:compred-circ},
and that by Proposition~\ref{p:lambdacirc} $\lambda(A)$ is 
the maximum cycle mean of any maximal strongly connected component of
$\digr(A)$.

Since $\lambda(A)=0$ implies $A=0$ by Corollary~\ref{l:zerocase},
we can assume $\lambda(A)=1$ without loss of generality. 
Since all nodes of $\digr(A)$ are critical, the transient of periodicity of $\{A,A^2,A^3,\ldots\}$ is the same
as the greatest transient of periodicity of any sequence of rows of these powers $\{A_{i\bullet},A^2_{i\bullet},A^3_{i\bullet},\ldots\}$ where $i$ is critical. However, these transients
are bounded by $(n-1)^2+1$ by~\cite{MNSS} Main Theorem~1. 
\kd

The following proposition gives more 
information on the critical graph and cyclicity of circulant matrices.

\begin{proposition}
\label{t:critcirc}
Let $A=\Circulant(a_0,\dots,a_{n-1})\neq 0$ and let $p_1,\ldots,p_s\in\{1,\ldots,n-1\}$
be the nonzero indices for which $a_{p_1}=\ldots=a_{p_s}=\lambda(A)$ (if such indices exist) and
such that $p_1>p_2>\ldots>p_s$. Then\\ 
\begin{itemize}
\item[(i)] $\crit(A)$ consists of $m=\gcd(n,p_1,\ldots,p_s)$ isomorphic strongly connected components.
Node set of the $i$th component, for $i\in\{1,\ldots,m\}$, is $\{i,i+m,\ldots,i+(n/m-1)m\}$.


\item[(ii)] $\per(A)$, equal to the cyclicity of each of these components, is $1$ if $a_0=\lambda(A)$
and 
\begin{equation}
\label{e:circcycl}
\begin{split}
&\per(A)=\gcd(\frac{n}{\gcd(n,p_1)},\frac{p_1-p_2}{\gcd(p_1,p_2)}, \frac{p_1-p_3}{\gcd(p_1,p_3)},\ldots,\frac{p_1-p_s}{\gcd(p_1,p_s)})\\ 
&=\gcd(\frac{n}{\gcd(n,p_1)},\frac{p_1-p_2}{\gcd(p_1,p_2)}, \frac{p_2-p_3}{\gcd(p_2,p_3)}, \ldots, \frac{p_{s-1}-p_s}{\gcd(p_{s-1},p_s)})\\
&=\gcd(\frac{n}{\gcd(n,p_1)},\frac{p_1-p_2}{\gcd(n,p_1,p_2)}, 
\frac{p_1-p_3}{\gcd(n,p_1,p_2,p_3)} \ldots, \frac{p_1-p_s}{\gcd(n,p_1,\ldots,p_s)})
\end{split}
\end{equation}
if $a_0\neq\lambda(A)$.
\end{itemize}
\end{proposition}

Parts of this statement can be found in~\cite{mp} Theorem 4.1 and Lemma 4.1.
Essentially, part (i) was proved in~\cite{mmp} Lemma 4.2 and Lemma 4.3, although
in the max-min algebra setting. The number $\gcd(n,p_1,\ldots,p_s)$ also appeared in~\cite{Tom-10} Theorem 4 as the "eigenspace dimension". The result of part (ii) relies on~\cite{Gav-99} (Theorems 3.1 and 3.3) 
where the cyclicity of threshold circulant graphs (see Definition~\ref{def:thresh}) 
was studied. 
We will give a complete proof of (i) and a reduction of (ii) to the results of~\cite{Gav-99} in
Subsection~\ref{ss:proof1}, for the reader's convenience.

Let us now describe the attraction cone of a circulant matrix as a solution set of 
a max-algebraic two-sided system of equations.
\begin{proposition}
\label{p:attrdesc-circ}
Let $A\in\Rpnn$ be a circulant matrix. Then
$$\Attr(A)=\{x\colon \lambda(A)A^{n^2}\otimes x=A^{n^2+1}\otimes x\}.$$ 
\end{proposition}
{\it Proof.} 
By Corollary~\ref{l:zerocase}, $\lambda(A)=0$ if and only if $A=0$, in which case
$\attr(A)=\Rpn$ and $A^{n^2}=A^{n^2+1}=0$, and the claim holds trivially.
Otherwise, by Corollary~\ref{c:compred-circ} $A$ is completely reducible and by Proposition~\ref{p:lambdacirc}
the maximal cycle mean of each component of $\digr(A)$ is the same. The claim then follows since
$A$ satisfies the conditions of Proposition~\ref{p:attrdesc} and since $n^2\geq T(A)$ by
Proposition~\ref{p:transient-circ}.
\kd

Let us examine the attraction cone of a $4\times 4$ circulant matrix.

\begin{example}
\label{ex:4x4}
{\rm Consider
\begin{equation*}
A=
\begin{pmatrix}
0 & 0 & 1 & t\\
t & 0 & 0 & 1\\
1 & t & 0 & 0\\
0 & 1 & t & 0
\end{pmatrix},
\end{equation*}
where $t\colon 0<t<1$. This is a circulant matrix, $\lambda(A)=1$, 
and $\crit(A)$ consists of two disjoint cycles:
$(1\; 3)$ and $(2\; 4)$. The cyclicity of $\crit(A)$ is thus equal to $2$ 
and so is the ultimate period of the max-algebraic matrix powers of $A$. 
Taking the max-algebraic powers of $A$ we obtain
\begin{equation*}
\begin{split}
A^{2} &=
\begin{pmatrix}
1 & t & t^2 & 0\\
0 & 1 & t & t^2\\
t^2 & 0 & 1 & t\\
t & t^2 & 0 & 1
\end{pmatrix}, \quad 
A^{2k}=
\begin{pmatrix}
1 & t & t^2 & t^3\\
t^3 & 1 & t & t^2\\
t^2 & t^3 & 1 & t\\
t & t^2 & t^3 & 1
\end{pmatrix}\quad \forall k\geq 2,\\
A^{2k-1}&=
\begin{pmatrix}
t^2 & t^3 & 1 & t\\
t & t^2 & t^3 & 1\\
1 & t & t^2 & t^3\\
t^3 & 1 & t & t^2
\end{pmatrix}\quad \forall k\geq 2.
\end{split}
\end{equation*} 
In particular, the periodicity transient is $T(A)=3$. 
By Proposition~\ref{p:attrdesc-circ} we have
$\attr(A)=\{x\colon A^{16}\otimes x=A^{17}\otimes x\}$, implying that the attraction cone
is precisely the set of vectors $x=(x_1\;x_2\; x_3\;x_4)$ that satisfy
\begin{equation}
\label{e:ex-nonreduced}
\begin{split}
x_1\oplus tx_2\oplus t^2 x_3\oplus t^3x_4&=t^2x_1\oplus t^3x_2\oplus x_3\oplus tx_4\\
tx_1\oplus t^2 x_2\oplus t^3x_3\oplus x_4&=t^3x_1\oplus x_2\oplus tx_3\oplus t^2 x_4.
\end{split}
\end{equation}
System~\eqref{e:ex-nonreduced} can be further reduced using the cancellation rule
$$a\oplus b= ta\oplus c\Leftrightarrow a\oplus b=c,$$
where $t<1$ and $a,b,c$ are arbitrary. Repeatedly applying this rule we obtain the 
system
\begin{equation}
\label{e:ex-reduced}
\begin{split}
x_1\oplus tx_2&=x_3\oplus tx_4\\
tx_1\oplus x_4&=x_2\oplus tx_3,
\end{split}
\end{equation}
equivalent to~\eqref{e:ex-nonreduced}.

Now observe that $x=[t\; 1\; t^2\; 1]$ satisfies this system of equations and 
belongs to the attraction cone. In particular, the ultimate period of $\{A^tx\}_{t\geq 1}$ is $1$,
however, $A\otimes x\neq x$ which shows that $\attr(A)$ is not the same as the (max-algebraic)
eigencone of $A$ in this case.}
\end{example}

The following theorem is one of the main results of the paper.
Its proof is postponed to Subsection~\ref{ss:proof2}.

\begin{theorem}
	\label{t:attrincl}
	Let $A,B\in\Rpnn$ be two circulant matrices such that $\lambda(A)=\lambda(B)$ and $A\leq B$. Then
	$\Attr(A)\subseteq\Attr(B)$.
\end{theorem}

Let us give two examples demonstrating this theorem. In the first example we have
two $0$-$1$ matrices, and in the second one we consider the matrix of Example~\ref{ex:4x4}
with two different values of $t$.

\begin{example}
\label{ex:01}
{\rm
Let us first consider a pair of $0$-$1$ matrices:
\begin{equation*}
A=
\begin{pmatrix}
0 & 1 & 0 & 0 & 0 & 0\\
0 & 0 & 1 & 0 & 0 & 0\\
0 & 0 & 0 & 1 & 0 & 0\\
0 & 0 & 0 & 0 & 1 & 0\\
0 & 0 & 0 & 0 & 0 & 1\\
1 & 0 & 0 & 0 & 0 & 0
\end{pmatrix},\quad
B=
\begin{pmatrix}
0 & 1 & 0 & 1 & 0 & 0\\
0 & 0 & 1 & 0 & 1 & 0\\
0 & 0 & 0 & 1 & 0 & 1\\
1 & 0 & 0 & 0 & 1 & 0\\
0 & 1 & 0 & 0 & 0 & 1\\
1 & 0 & 1 & 0 & 0 & 0
\end{pmatrix}.
\end{equation*}
 Observe that the sequence $\{A^t\}_{t\geq 1}$ is periodic from the very beginning. The system 
$A^{36}\otimes x=A^{37}\otimes x$, being the same as $A\otimes x=A^2\otimes x$,
reduces to $x_1=x_2=x_3=x_4=x_5=x_6$.\\
The sequence $\{B^t\}_{t\geq 1}$ becomes periodic from $T(B)=2$. More precisely, we have
\begin{equation*}
B^{2k+1}=
\begin{pmatrix}
0 & 1 & 0 & 1 & 0 & 1\\
1 & 0 & 1 & 0 & 1 & 0\\
0 & 1 & 0 & 1 & 0 & 1\\
1 & 0 & 1 & 0 & 1 & 0\\
0 & 1 & 0 & 1 & 0 & 1\\
1 & 0 & 1 & 0 & 1 & 0
\end{pmatrix},\quad
B^{2k}=
\begin{pmatrix}
1 & 0 & 1 & 0 & 1 & 0\\
0 & 1 & 0 & 1 & 0 & 1\\
1 & 0 & 1 & 0 & 1 & 0\\
0 & 1 & 0 & 1 & 0 & 1\\
1 & 0 & 1 & 0 & 1 & 0\\
0 & 1 & 0 & 1 & 0 & 1
\end{pmatrix},\quad k\geq 1,
\end{equation*}
and the system $B^{36}\otimes x=B^{37}\otimes x$ reduces to 
$x_1\oplus x_3\oplus x_5=x_2\oplus x_4\oplus x_6$, thus $\attr(A)\subseteq\attr(B)$. 
\kd
}
\end{example}

\if{
$\crit(A)$ has $6$ cyclic classes, each of them consisting just of 
one node, and the attraction cone $\attr(A,1)$ consists of vectors $x$
satisfying $x_1=x_2=x_3=x_4=x_5=x_6$.\\
$\crit(B)$ has cyclicity $2$ and consists of two cyclic classes 
$\{1,3,5\}$ and $\{2,4,6\}$, and the attraction cone $\attr(B,1)$ consists of vectors $x$
satisfying $x_1\oplus x_3\oplus x_5=x_2\oplus x_4\oplus x_6$. 
Obviously, each vector whose components are equal to each other satisfies this system,
thus $\attr(A)\subseteq\attr(B)$.
}\fi

\begin{example}
\label{ex:4x42t}
{\rm
Take
\begin{equation*}
A= 
\begin{pmatrix}
0 & 0 & 1 & t_1\\
t_1 & 0 & 0  & 1\\
1 & t_1 & 0 & 0\\
0 & 1 & t_1 & 0
\end{pmatrix},\quad
B=
\begin{pmatrix}
0 & 0 & 1 & t_2\\
t_2 & 0 & 0  & 1\\
1 & t_2 & 0 & 0\\
0 & 1 & t_2 & 0
\end{pmatrix}
\end{equation*}
with $0<t_1<t_2<1$. Then $\attr(A)$ is the set of all $x$ satisfying~\eqref{e:ex-reduced}
with $t=t_1$, which is
\begin{equation}
\label{e:ex-reduced1}
\begin{split}
x_1\oplus t_1x_2&=x_3\oplus t_1x_4\\
t_1x_1\oplus x_4&=x_2\oplus t_1x_3,
\end{split}
\end{equation}
and $\attr(B)$ is the set of all $x$ satisfying
\begin{equation}
\label{e:ex-reduced2}
\begin{split}
x_1\oplus t_2x_2&=x_3\oplus t_2x_4\\
t_2x_1\oplus x_4&=x_2\oplus t_2x_3,
\end{split}
\end{equation}
We next show that $\attr(A)\subseteq\attr(B)$ in this example, by considering 
various special cases.

Suppose first that we have $t_1x_2=t_1x_4\geq x_1\oplus x_3$ in the first equation of~\eqref{e:ex-reduced1}.
This implies $x_2=x_4\geq t_2(x_1\oplus x_3)\geq t_1(x_1\oplus x_3)$ and $t_2x_2=t_2x_4\geq (x_1\oplus x_3)$.
This shows that in this case $x$ belongs to both $\attr(A,1)$ and $\attr(B,1)$. The case when 
$t_1x_1=t_1x_3\geq x_2\oplus x_4$ in the second equation of~\eqref{e:ex-reduced1} is treated similarly.

Suppose now that $x\in\Attr(A)$ and $t_1x_2=x_3\geq x_1\oplus t_1 x_4$. As we cannot have
$t_1x_1=x_2$ and $x_4=t_1x_3$ in the second equation of~\eqref{e:ex-reduced1}, assume that 
$x_2=x_4\geq t_1(x_1\oplus x_3)$. But this implies $t_1x_2=t_1x_4$, and as $t_1x_2$ is the maximum in the first equation, this returns us to the case which we considered first, where $x\in\Attr(B)$. We also note three other similar cases that are treated in the same way.

The remaining case when $x\in \Attr(A)$, $x_1=x_3\geq t_1(x_2\oplus x_4)$ and $x_2=x_4\geq t_1(x_1\oplus x_3)$ is impossible when $t_1<1$.}\kd 
\end{example}

\section{Interval robustness of circulant matrices}

\label{s:interval}

In this  section we characterize the six types of interval robustness of Definition \ref{def:introb} for interval circulant matrix $\IntCirculant(\bfa_0,\ldots,\bfa_{n-1})$ and interval $\bfX=\times_{i=1}^n \bfX_i$ where $\bfX_i$ and $\bfa_i$ are intervals independently taking one of the following four forms:
$$[\undx_i,\ovex_i],\ (\undx_i,\ovex_i),\ (\undx_i,\ovex_i],\ [\undx_i,\ovex_i)$$
and 
$$[\unda_j,\ovea_j],\ (\unda_j,\ovea_j),\ (\unda_j,\ovea_j],\ [\unda_j,\ovea_j)$$
for $\undx_i,\ovex_i\in\Rp$ and $i\in N$, and $\unda_j,\ovea_j\in\Rp$ and $j\in N_0$, respectively.


\subsection{Universal and possible  $\bfX-$robustness} 

\label{ss:univposs}

Let us introduce the following notation.

\begin{definition}[Matrices $A^{(k)}$ and vectors $x^{(k)}$] 
\label{def:Akxk}
{\rm For a given index  $k\in \Np$ denote $$A^{(k)}=\Circulant(\unda_0,\unda_1,\dots,\unda_{k-1},\ovea_k,\unda_{k+1},\dots,\unda_{n-1}),$$
and 
$$
x^{(k)}=(\undx_1,\undx_2,\dots,\undx_{k-1},\ovex_k,\undx_{k+1},\dots,\undx_{n})
$$
}
\end{definition}

\if{
\begin{definition}[Intervals $\mbf{\check{X}}$ and $\mbf{\check{X}}^A$]
\label{def:AvXv}
{\rm Denote:}
$$\mbf{\check{X}}_i=\{s\colon s\cdot \overline{x}_i\in\mbf{X}_i\},\quad
\mbf{\check{X}}=\times_{i=1}^n \mbf{\check{X}}_i,$$
$$\mbf{\check{X}}^A_i=\{s\colon s\cdot\overline{a}_i\in\mbf{a}_i\},\quad
\mbf{\check{X}}^A=\times_{i=1}^n \mbf{\check{X}}^A_i.$$
\end{definition}
}\fi

\if{
{\color{red}(For the symmetry of  notations I suggest to denote $\mbf{\check{a}}_i$ instead of $\mbf{\check{X}}^A_i$)}
}\fi


\if{
\begin{lemma}\label{t:L1}
	Let $\mbf{X}\subseteq\Rpn$ be an interval vector. 
	Then $x\in\mbf{X}$ if and only if $x=\bigoplus_{k=1}^n \beta_k x^{(k)}$
	with $\beta\in\mbf{\check{X}}$. In particular,
	\begin{equation}
	\label{e:xgen}
	x=\bigoplus_{k=1}^n \frac{x_k}{\overline{x}_k} x^{(k)}.
	\end{equation}
\end{lemma}	
{\it Proof.} Let us first show that~\eqref{e:xgen} holds, meaning that it holds 
for each coordinate $x_i$.
Observe that when $k\neq i$ we have that  $x_k/\ovex_k\leq 1$ implies 
$(x_k/\ovex_k) x_i^{(k)}\leq\undx_i$, and when $k=i$ we obtain
$(x_i/\ovex_i) x_i^{(i)}=x_i$. Since $\undx_i\leq x_i$, we obtain that indeed
\begin{equation*}
\bigoplus_{k=1}^n \frac{x_k}{\overline{x}_k} x^{(k)}_i=(x_i/\ovex_i) x_i^{(i)}=x_i,
\end{equation*}
so~\eqref{e:xgen} holds. 

"Only if": When $x\in\mbf{X}$ we can choose $\beta_k=\frac{x_k}{\overline{x}_k}\in\mbf{\check{X}_k}$ for all $k$. 

"If": Observe that
$\beta_i x^{(i)}_i=\beta_i\ovex_i\in \mbf{X}_i$ by the definition of $\mbf{\check{X}_i}$,
while $\beta_k x^{(k)}_i=\beta_k \undx_i\leq \undx_i\leq \beta_i x^{(i)}_i$ for $k\neq i$. Therefore
$x_i=\bigoplus_{k=1}^n \beta_k x^{(k)}_i\in\mbf{X}_i$ for each $i$ and $x\in\mbf{X}$.
\kd	
}\fi

\if{	
\begin{lemma}\label{t:L2}
	Let  $\IntCirculant({\mbf{a_0},\mbf{a_1},\ldots,\mbf{a_{n-1}}})$ be an interval circulant matrix. 
	Then $A\in\IntCirculant(\bfa_0,\dots,\bfa_{n-1})$ if and only if $A=\bigoplus_{k=0}^{n-1} \gamma_k A^{(k)}$
	with $\gamma\in\mbf{\check{X}}^A$. In particular,
	\begin{equation}
	\label{e:Agen}
	A=\bigoplus_{k=1}^n \frac{a_k}{\overline{a}_k} A^{(k)}.
	\end{equation}
\end{lemma}		 
{\it Proof.} 
Similar to the proof of Lemma~\ref{t:L1}.
\kd	
}\fi

The following lemma explains the use of vectors $x^{(k)}$. 

\begin{lemma}
	\label{l:generators}
	Let $\bfX\subseteq\Rpn$ be an interval and let $A\in\Rpnn$. Then $\mbf{X}\subseteq\attr(A)$ 
	if and only if $x^{(i)}\in\Attr(A)$ for each $i\in N.$
\end{lemma}
{\it Proof.}
	Observe first that since the cone $\attr(A)$ is a closed set by Corollary~\ref{c:attrclosed}, the inclusion 
	$\mbf{X}\subseteq\Attr(A)$ is equivalent to $\operatorname{cl}(\mbf{X})\subseteq\Attr(A)$, where $\operatorname{cl}$ is a Euclidean closure. 
Since $x^{(i)}\in\operatorname{cl}(\mbf{X})$ for all $i\in N$ (as vertices of 
the box	$\operatorname{cl}(\mbf{X})$), it follows that the condition is necessary. Let us show 
that this condition is also sufficient.
For this we will show that
\begin{equation}
	\label{e:xgen}
	x=\bigoplus_{k=1}^n \frac{x_k}{\overline{x}_k} x^{(k)}.
	\end{equation}	
Indeed, observe that when $k\neq i$ we have that  $x_k/\ovex_k\leq 1$ implies 
$(x_k/\ovex_k) x_i^{(k)}\leq\undx_i$, and when $k=i$ we obtain
$(x_i/\ovex_i) x_i^{(i)}=x_i$. Since $\undx_i\leq x_i$, we obtain that 
\begin{equation*}
\bigoplus_{k=1}^n \frac{x_k}{\overline{x}_k} x^{(k)}_i=(x_i/\ovex_i) x_i^{(i)}=x_i,
\end{equation*}
for all $i$, so~\eqref{e:xgen} holds. Thus $x$ can be expressed as a max-linear combination
of $x^{(k)}$ for $k\in N$ and $x\in\attr(A)$ since $\attr(A)$ is a max-cone
(Definition~\ref{def:maxcone}).
\kd

\begin{definition}[Matrix $\hat{A}$]
\label{def:hatA}
For  $\unda=\max\limits_{k\in \Np}\unda_k$   define  
$$\Hat{A}=\Circulant(\hat{a}_0,\hat{a}_1,\dots,\hat{a}_{n-1}),$$ 
where
\begin{equation*}\label{E_PR}
\hat{a}_i=\min\{\unda,\ovea_i\},\text{   for each   } i\in \Np.
\end{equation*}
\end{definition}

 Let us characterize the cases when $\hat{A}=0$ and when 
$\hat{A}\in\IntCirculant(\bfa_0,\dots,\bfa_{n-1})$.
\begin{proposition}
\label{p:Ahat} Let $\IntCirculant(\bfa_0,\dots,\bfa_{n-1})$ be given. Then
\begin{itemize}
\item[{\rm(i)}] $\hat{A}=0\Leftrightarrow\underline{a}=0\Leftrightarrow \underline{A}=0\Leftrightarrow \lambda(\underline{A})=0.$
\item[{\rm(ii)}] If $\forall i:\bfa_i=[\underline{a_i},\overline{a_i}]$, then $\hat{A}\in \IntCirculant(\bfa_0,\dots,\bfa_{n-1}).$
\item[{\rm(iii)}] $\hat{A}\in \IntCirculant(\bfa_0,\dots,\bfa_{n-1}) \Leftrightarrow \forall i: (\unda\geq\ovea_i\Rightarrow \ovea_i\in \bfa_i)\&(\unda\leq\ovea_i\Rightarrow \unda\in \bfa_i).$
\end{itemize}
\end{proposition}
{\it Proof.}
(i): Let us show that $\hat{A}=0\Leftrightarrow\underline{a}=0$. By Definition~\ref{def:hatA} it is 
immediate that $\underline{a}=0$ implies $\hat{A}=0$. Next, assume that $\hat{A}=0$. 
Then $\overline{a}_i=0$ for all $i$, which implies $\underline{a}_i=0$ for all $i$, hence $\underline{a}=0$. 
The equivalence $\underline{a}=0\Leftrightarrow \underline{A}=0$ is obvious, and 
$\underline{A}=0\Leftrightarrow \lambda(\underline{A})=0$ follows from Corollary~\ref{l:zerocase}.\\
(ii) and (iii): Straightforward.
\kd

Matrices $\hat{A}$ and $A^{(k)}$ for $k=0,\ldots,n-1$ have the following useful properties.

\begin{lemma}
\label{l:Ahat}
If $\hat{A}\neq 0$, then $(\forall A\in\IntCirculant(\bfa_0,\dots,\bfa_{n-1}))[(A/\lambda(A))\leq (\hat{A}/\lambda(\hat{A})]$.
\end{lemma}
{\it Proof.}
Observe that $\hat{A}\neq 0$ implies that $A=0$ does not belong to the interval matrix $\IntCirculant(\bfa_0,\dots,\bfa_{n-1})$.
Recalling that $\hat{a}_i=\min(\overline{a}_i,\underline{a})$ for all $i$ 
	we see that $\hat{a}_i\leq\underline{a}$ for all $i$ and that 
	$\hat{a}_k=\underline{a}$ for $k$ such that $\underline{a}_k=\underline{a}$.
	Hence $\lambda(\hat{A})=\underline{a}$ by Proposition~\ref{p:lambdacirc}. Showing $(A/\lambda(A))\leq (\hat{A}/\lambda(\hat{A})$ means showing 
	\begin{equation}
	\label{e:toprove}
	a_i/\max_k a_k \leq \min(\overline{a}_i,\underline{a})/\underline{a}\quad\forall i.
	\end{equation}
	To prove~\eqref{e:toprove} we observe that it follows from the inequality
	\begin{equation}
	\label{e:toprove1}
	a_i\cdot \underline{a}\leq\max_j a_j\cdot\min(\overline{a}_i,\underline{a})\quad\forall i,
	\end{equation}
	which is 
	\begin{equation}
	\label{e:case1}
	a_i\cdot\underline{a}\leq\max_j a_j\cdot\underline{a}
	\end{equation}
	when $\min(\overline{a}_i,\underline{a})=\underline{a}$, and 
	\begin{equation}
	\label{e:case2}
	a_i\cdot \max_i{\underline{a}_i}\leq\overline{a}_i\cdot \max_j a_j
	\end{equation}
  when $\min(\overline{a}_i,\underline{a})=\overline{a}_i$.
	Both~\eqref{e:case1} and~\eqref{e:case2} are obvious.
	This shows~\eqref{e:toprove1} and hence~\eqref{e:toprove} and 
	$(A/\lambda(A))\leq (\hat{A}/\lambda(\hat{A}))$. 
\kd

\begin{lemma}
\label{l:Akleq}
For any nonzero $A\in\IntCirculant(\bfa_0,\dots,\bfa_{n-1})$
there exists $A^{(k)}\neq 0$ for some $k\in N_0$ such that 
$[(A^{(k)}/\lambda(A^{(k)}))\leq (A/\lambda(A))]$.
\end{lemma}
{\it Proof.}
Let $A=\Circulant(a_0,\dots,a_{n-1})$ and let $k$ be such that 
$a_k=\max_{j\in N} a_j$. Consider $A^{(k)}$. Since $ \overline{a}_k\geq a_k>0$ but the rest of the components defining $A^{(k)}$ are $\underline{a}_i\leq a_i$ for $i\neq k$, we have $\lambda(A^{(k)})=\overline{a}_k$ and 
$(A^{(k)}/\lambda(A^{(k)}))\leq (A/\lambda(A))$.
\kd

We now characterize possibly $\mbf{X}$-robust and
universally $\mbf{X}$-robust interval circulant matrices. 

\begin{theorem}
	\label{t:posrob}
	Let $\mbf{X}\subseteq\Rpn$ be an interval, and let  
	$\IntCirculant(\bfa_0,\dots,\bfa_{n-1})\subseteq\Rpnn$ be an interval circulant matrix 
	containing $\hat{A}$. Then $\IntCirculant(\bfa_0,\dots,\bfa_{n-1})$ is possibly $\mbf{X}$-robust
	if and only if we have $x^{(i)}\in\attr(\hat{A})$ for all $i\in N$.
\end{theorem}
{\it Proof.}
	We need to show that there exists $A\in\IntCirculant(\bfa_0,\dots,\bfa_{n-1})$ such that
	$\mbf{X}\subseteq\attr(A)$ 
	if and only if $x^{(i)}\in\attr(\hat{A})$ for all $i\in N$.
	If $\hat{A}=0$ then $\attr(\hat{A})=\Rpn$ and the claim is obvious. Next we suppose that
	$\hat{A}\neq 0$ which implies $\lambda(\hat{A})\neq 0$ by Corollary~\ref{l:zerocase}.
	By Proposition~\ref{p:Ahat} part (i), $\IntCirculant(\bfa_0,\dots,\bfa_{n-1})$ contains only nonzero 
	matrices in this case.
	
	{\bf ``If'':} By Lemma~\ref{l:generators}, the condition implies that $\mbf{X}\subseteq\attr(\hat{A})$. The claim then follows since $\hat{A}\in\IntCirculant(\bfa_0,\dots,\bfa_{n-1})$.

	{\bf ``Only if'':} Let $A\in\IntCirculant(\bfa_0,\dots,\bfa_{n-1})$ be such that 
	$\mbf{X}\subseteq\attr(A)$. By Lemma~\ref{l:Ahat} we have $(A/\lambda(A))\leq (\hat{A}/\lambda(\hat{A}),$
	and Theorem~\ref{t:attrincl} yields that~$x\in\Attr{\hat{A}}$.  
	As $x\in\attr{\hat{A}}$ for all $x\in\mbf{X}$, 
	the claim then follows from Lemma~\ref{l:generators}.
	\kd

\begin{corollary}\label{cor-pos-robust} 
	Let $x\in \Rpn$ and let $\IntCirculant(\bfa_0,\dots,\bfa_{n-1})\subseteq\Rpnn$ be an interval
	circulant matrix containing 
	$\hat{A}$. Then  
	$(\exists A\in\IntCirculant(\bfa_0,\dots,\bfa_{n-1}))[x\in\attr(A)]$ if and only if
	$x\in\attr(\hat{A})$.
\end{corollary}
{\it Proof.} Take $\mbf{X}=\{x\}$ then the possible $\mbf{X}$-robustness means existence of
$A\in \IntCirculant(\bfa_0,\dots,\bfa_{n-1})$ such that $x\in\attr(A)$ and $x^{(i)}=x$ for all $i\in N$.
The claim then follows from Theorem~\ref{t:posrob}.
\kd

\begin{theorem}
	\label{t:unirob}
Let $\mbf{X}\subseteq\Rpn$ be an interval, and let $\IntCirculant(\bfa_0,\dots,\bfa_{n-1})\subseteq\Rpnn$
	be an interval circulant matrix. Then $\IntCirculant(\bfa_0,\dots,\bfa_{n-1})$ is universally $\mbf{X}$-robust if and only if  $x^{(j)}\in\Attr(A^{(i)})$ for all $i\in\Np$ and $j\in N$.
\end{theorem}
{\it Proof.}
	We need to show that $\mbf{X}\subseteq\attr(A)$ for all $A\in\IntCirculant(\bfa_0,\dots,\bfa_{n-1})$ 
	if and only if $x^{(j)}\in\Attr(A^{(i)})$ for all $i\in N_0$ and $j\in N$. 
	
 {\bf ``If'':} 
	Let $x^{(j)}\in\attr(A^{(i)})$ hold for all $i\in N_0$ and $j\in N$.
 Take $A\in\IntCirculant(\bfa_0,\dots,\bfa_{n-1})$. If $A=0$ then $x^{(j)}\in\attr(A)=\Rpn$.
Otherwise, by Lemma~\ref{l:Akleq} there exists $k\in N_0$ such
	that $A^{(k)}\neq 0$ and $(A^{(k)}/\lambda(A^{(k)}))\leq (A/\lambda(A))$.
	Applying Theorem~\ref{t:attrincl} to $(A^{(k)}/\lambda(A^{(k)}))$ and $(A/\lambda(A))$ we obtain $x^{(j)}\in\attr(A)$ for all nonzero $x^{(j)}$, hence $\mbf{X}\subseteq\attr(A)$.

	{\bf ``Only if'':} Take a sequence $\{A_s\}_{s\geq 1}\subseteq \IntCirculant(\bfa_0,\dots,\bfa_{n-1})$ such that
	$\lim_{s\to\infty} A_s=A^{(k)}$, and take any $x\in\mbf{X}$. Since $x\in\attr(A_s)$ for all $s$, by Proposition~\ref{p:attrdesc-circ} we have
	$\lambda(A_s)A_s^{n^2}\otimes x=A_s^{n^2+1}\otimes x$ for all $s$, and by the continuity of the arithmetic
	operations of max-algebra we obtain
	$\lambda(A^{(k)})(A^{(k)})^{n^2}\otimes x=(A^{(k)})^{n^2+1}\otimes x$.
	As $x\in\attr(A^{(k)})$ for all $x\in\mbf{X}$, 
	the claim then follows from Lemma~\ref{l:generators}.
\kd

\begin{corollary} \label{cor_uni_robust}
	Let $x\in\Rpn$, and let $\IntCirculant(\bfa_0,\dots,\bfa_{n-1})\subseteq\Rpnn$
	be an interval circulant matrix. . Then
	$(\forall A\in \IntCirculant(\bfa_0,\dots,\bfa_{n-1})) \, [x\in\attr(A) \,]$   
	if and only if   $x\in\attr(A^{(k)})$ for each $k\in \Np$.
\end{corollary}

{\it Proof.} 
Take $\mbf{X}=\{x\}$ then the universal $\mbf{X}$-robustness means that $x\in\attr(A)$ for all
$A\in \IntCirculant(\bfa_0,\dots,\bfa_{n-1})$.
The claim then follows from Theorem~\ref{t:unirob}.
\kd

\subsection{Tolerance and weak tolerance $\bfX-$robustness} 

\label{ss:tolweak}
 
\begin{theorem}\label{TXRCIR}
	Let $\bfX\subseteq\Rpn$ be a closed interval, and let $\IntCirculant(\bfa_0,\dots,\bfa_{n-1})\subseteq\Rpnn$
	be an interval circulant matrix.
	Then $\IntCirculant(\bfa_0,\dots,\bfa_{n-1})$ is  
	tolerance $\bfX-$robust if and only if $(\forall k\in N_0)[(\attr(A^{(k)})\cap\mbf{X})\neq\emptyset]$.
\end{theorem}
 {\it Proof.}
{\bf ``If'':} 
Take $A\in\IntCirculant(\bfa_0,\dots,\bfa_{n-1})$. If $A=0$ then $\attr(A)=\Rpn$, hence $\attr(A)\cap\mbf{X}\neq\emptyset$. Otherwise, for each $i\in N_0$ take $y^{(i)}\in(\bfX\cap \attr(A^{(i)})]$, 
By Lemma~\ref{l:Akleq} there exists $k\in N_0$ 
with $(A^{(k)}/\lambda(A^{(k)}))\leq (A/\lambda(A))$. Applying Theorem~\ref{t:attrincl} to $(A^{(k)}/\lambda(A^{(k)}))$ and $(A/\lambda(A))$ we obtain $y^{(k)}\in\attr(A)$, hence the implication.

{\bf ``Only if'':}
For any $k\in N_0$ take a sequence $\{A_s\}_{s\geq 1}\subseteq \IntCirculant(\bfa_0,\dots,\bfa_{n-1})$ such that
	$\lim_{s\to\infty} A_s=A^{(k)}$. 
	For each of these matrices there exists $x^s\in\mbf{X}$ such that
	$x^s\in\attr(A_s)$. Then by Proposition~\ref{p:attrdesc-circ} we have $\lambda(A_s)A_s^{n^2}\otimes x^s=A_s^{n^2+1}\otimes x^s$ for all $s$.
	Since $\mbf{X}$ is compact, we can assume that $\lim_{s\to\infty} x^s$ exists and denote it by $y^{(k)}$. Then we obtain that
	by the continuity of arithmetic operations of max-algebra $\lambda(A^{(k)})(A^{(k)})^{n^2}\otimes y^{(k)}=(A^{(k)})^{n^2+1}\otimes y^{(k)}$. Hence
	$y^{(k)}\in\Attr(A^{(k)})$.
\kd

\begin{corollary}
\label{CXRCIR}
Under the conditions of Theorem~\ref{TXRCIR}, $\IntCirculant(\bfa_0,\dots,\bfa_{n-1})$ is tolerance
$\bfX-$robust if and only if all systems
\begin{equation}
\label{EXRCIR}
\lambda(A^{(k)})(A^{(k)})^{n^2}\otimes y=(A^{(k)})^{n^2+1}\otimes y,\quad y\in\mbf{X},
\end{equation}
with $k\in N_0$ such that $A^{(k)}\neq 0$ are solvable.
\end{corollary}

We now characterize the weak tolerance robust matrices.

\if{
\begin{definition}[Matrix $\Akr(\bfX,r)$]
\label{def:akr}
Let $\IntCirculant(\bfa_0,\dots,\bfa_{n-1})$ be an interval circulant matrix such that  $\underline{A}\neq 0$ and   $\hat{A}\in\IntCirculant(\bfa_0,\dots,\bfa_{n-1})$, and let $\bfX$ be an interval. 
Define  matrix $\Akr(\bfX,r)\in \Rpnn:$ 
$$ \Akr(\bfX,r) =\left(\begin{array}{ccc}
\left(\frac{\hat{A}}{\lambda(\hat{A})}\right)^r\otimes x^{(1)},\dots,
\left(\frac{\hat{A}}{\lambda(\hat{A})}\right)^r\otimes x^{(n)} \\
\end{array}\right).$$
\end{definition}
}\fi

\begin{theorem}\label{T_weak_tolerance_robust} 
Let $\bfX\subseteq\Rpn$ be an interval and let 
$\IntCirculant(\bfa_0,\dots,\bfa_{n-1})\subseteq\Rpnn$ be an interval circulant matrix 
containing $\hat{A}$.
 Then $\IntCirculant(\bfa_0,\dots,\bfa_{n-1})$ is weakly tolerance $\bfX-$robust if and only if 
	$ \lambda(\hat{A})(\hat{A})^{n^2}\otimes x=
	  (\hat{A})^{n^2+1} \otimes x$
	is solvable with $x\in\mbf{X}$. 
\end{theorem}
{\it Proof.} 
By Corollary \ref{cor-pos-robust}, $x\in\mbf{X}$ and $A\in\IntCirculant(\bfa_0,\dots,\bfa_{n-1})$ such 
that $x\in\attr(A)$ exist if and only if $x\in\attr(\hat{A})$ for some $x\in\mbf{X}$. This,
by Proposition~\ref{p:attrdesc-circ}, is equivalent to 
$ \lambda(\hat{A})(\hat{A})^{n^2}\otimes x=
	  (\hat{A})^{n^2+1} \otimes x$
	being solvable with $x\in\mbf{X}$.
\kd

\subsection{Possible and tolerance $\IntCirculant(\bfa_0,\dots,\bfa_{n-1})-$robustness}
\label{ss:posstol}
\if{
Let us first introduce the following notation:
\begin{definition}[Matrix $\Dkr(\IntCirculant, r)$]
\label{def:Dkr} 
Let $\IntCirculant(\bfa_0,\dots,\bfa_{n-1})$ and  $\bfX$ be given. Define
\begin{displaymath}\hspace{-.8cm}
\Dkr(\IntCirculant,r) =  \left(\begin{array}{c}
\left(\frac{A^{(0)}}{\lambda(A^{(0)})}\right)^r\\
\left(\frac{A^{(1)}}{\lambda(A^{(1)})}\right)^r\\
\vdots\\
\left(\frac{A^{(n-1)}}{\lambda(A^{(n-1)})}\right)^r
\end{array}\right).\end{displaymath}
\end{definition}
}\fi

We now characterize the remaining two types of robustness.

\begin{theorem}\label{T_possiblyA_robust} 
	Let $\mbf{X}\subseteq\Rpn$ be an interval, and let $\IntCirculant(\bfa_0,\dots,\bfa_{n-1})\subseteq\Rpnn$
	be an interval circulant matrix. Then $\bfX$  is possibly $\IntCirculant(\bfa_0,\dots,\bfa_{n-1})-$robust  if and only if there exists $x\in\mbf{X}$ that satisfies
	$\lambda(A^{(i)})(A^{(i)})^{n^2}\otimes x=
		(A^{(i)})^{n^2+1}\otimes x$ for all $i\in N_0$
		such that $A^{(i)}\neq 0$.
\end{theorem}

	{\it Proof.} 
	By Corollary~\ref{cor_uni_robust}, $x\in\mbf{X}$ belongs to $\attr(A)$ for all $A\in\IntCirculant(\bfa_0\dots,\bfa_{n-1})$ if 
	and only if it belongs to $\attr(A^{(i)})$ for all $i\in N_0$ with $A^{(i)}\neq 0$. 
	By Proposition~\ref{p:attrdesc-circ} this is equivalent to $x$ satisfying $\lambda(A^{(i)})(A^{(i)})^{n^2}\otimes x=(A^{(i)})^{n^2+1}\otimes x$ for all such $i$.
	\kd

\begin{theorem}\label{T_toleranceA_robust} 
Let $\mbf{X}\subseteq\Rpn$ be an interval, and let $\IntCirculant(\bfa_0,\dots,\bfa_{n-1})\subseteq\Rpnn$
	be an interval circulant matrix containing 
$\hat{A}$. Then interval vector 
$\bfX$  is tolerance $\IntCirculant(\bfa_0,\dots,\bfa_{n-1})-$robust  if and only if $\IntCirculant(\bfa_0,\dots,\bfa_{n-1})$  is possibly $\bfX-$robust.
\end{theorem}
{\it Proof.} Suppose that
$\bfX$  is tolerance $\IntCirculant(\bfa_0,\dots,\bfa_{n-1})-$robust, then we have the following  
$$(\forall x\in \bfX) (\exists A\in \IntCirculant(\bfa_0,\dots,\bfa_{n-1})) [\, x\in\attr(A) \,]\stackrel{\mathrm{Cor. \ref{cor-pos-robust}}\vrule depth1mm width0mm height0mm}{\Longleftrightarrow}
 (\forall x\in \bfX) [\, x\in\attr(\hat A) \,]$$ 
$$\hspace{-.4cm} \Rightarrow 
(\exists A\in \IntCirculant(\bfa_0,\dots,\bfa_{n-1})) (\forall x\in \bfX) [\, x\in\attr(A) \,], $$ 
and hence we have that $\IntCirculant(\bfa_0,\dots,\bfa_{n-1})$  is possibly $\bfX-$robust.

The converse implication is trivial.
\kd

\subsection{Computational complexity}

We close the section with a couple of remarks on the computational complexity of 
the different types of interval robustness.

\begin{remark}
\label{r:complexity1}
{\rm By Theorems~\ref{t:posrob} and~\ref{T_toleranceA_robust}
the verification of whether
\begin{itemize}
\item[(i)] $\IntCirculant(\bfa_0,\dots,\bfa_{n-1})$ is possibly $\bfX$-robust,
\item[(ii)] $\IntCirculant(\bfa_0,\dots,\bfa_{n-1})$ is universally $\bfX$-robust,
\item[(iii)] $\bfX$ is tolerance $\IntCirculant(\bfa_0,\dots,\bfa_{n-1})$-robust
\end{itemize}
reduces, under some assumptions, to the verification whether some vectors satisfy some
two-sided max-linear systems with $n^2$ and $n^2+1$ powers of some matrices. Hence these types of 
robustness are of polynomial complexity.}
\end{remark}
\if{
Cases (i) and (iii) are equivalent and require the verification of $n$ such systems 
taking $O(n^3)$ operations. Computation of powers can be done by means of repeated squaring and 
takes $O(n^3\log n)$ operations~\cite{Ser-09}.
TODO: is it better for circulants??
Thus the overall complexity of (ii) and (iii)
(under the assumption that $\hat{A}\in \IntCirculant(\bfa_0,\dots,\bfa_{n-1})$) is $O(n^3\log n)$.

In Case (ii) we need to verify $n^2$ max-linear two-sided systems with powers of 
$n$ different matrices which raises the complexity to $O(n^4\log n)$.\\
TODO: can we do better with circulants?
}\fi

\begin{remark}
\label{r:complexity2}
{\rm By Corollary~\ref{CXRCIR}, Theorem~\ref{T_possiblyA_robust} and 
Theorem~\ref{T_toleranceA_robust}, verifying whether
\begin{itemize}
\item[(i)] $\IntCirculant(\bfa_0,\dots,\bfa_{n-1})$ is tolerance $\bfX$-robust,
\item[(ii)] $\IntCirculant(\bfa_0,\dots,\bfa_{n-1})$ is weakly tolerance $\bfX$-robust,
\item[(iii)] $\bfX$ is possibly $\IntCirculant(\bfa_0,\dots,\bfa_{n-1})$-robust
\end{itemize}
reduces, under some assumptions, to verifying the non-emptyness of solution set 
of some system of max-affine inequalities, where some of the inequalities
(among those defining $\bfX$) can be strict.  This problem was generally shown to be
polynomially equivalent to solving a mean-payoff game~\cite{All+}, for which efficient
pseudopolynomial algorithms exist, but existence of a polynomial algorithm has been a long-standing open question.}
\end{remark}

\section{Proofs of Proposition~\ref{t:critcirc} and Theorem~\ref{t:attrincl}}
\label{s:proofs}

\subsection{Cyclicity of circulants: Proof of Proposition~\ref{t:critcirc}}

\label{ss:proof1}

Let us start with the following elementary but useful statement.

\begin{lemma}
\label{l:NT}
Let $p_1,\ldots,p_s,n\in\Nat$ (the set of natural numbers). Then the equation
\begin{equation}
\label{e:p1ps}
p_1x_1+\ldots+p_sx_s\equiv m(\modd\; n)
\end{equation}
has a solution $(x_1,\ldots,x_s)\in(\Nat\cup\{0\})^s$ if and only if
$m$ is a multiple of $\gcd(p_1,\ldots,p_s,n)$.
\end{lemma}
{\it Proof.}
{\bf ``Only if"}: Observe that $p_1x_1+\ldots+p_sx_s$ 
and $n$ are always multiples of $\gcd(p_1,\ldots,p_s,n)$, and if~\eqref{e:p1ps} holds then so is $m$ as
well.

{\bf ``If"}: The claim is well known for $s=1$ (elementary number theory). 
The same fact also implies existence of $x_s\in\Nat\cup\{0\}$ such that
\begin{equation}
\label{e:ps}
p_sx_s\equiv m(\modd\;\gcd(n,p_1,\ldots,p_{s-1})).
\end{equation}
We now prove the claim by induction assuming that it holds for $s-1$. 
Observe that~\eqref{e:ps} implies that there exists also $k\in\Nat\cup\{0\}$ such that
\begin{equation}
\label{e:psk}
p_sx_s+k\gcd(n,p_1,\ldots,p_{s-1})\equiv m(\modd\;n).
\end{equation}
But by induction there exist $x_1\in\Nat\cup\{0\},\ldots,x_{s-1}\in\Nat\cup\{0\}$ such that
\begin{equation}
\label{e:p1ps-1}
p_1 x_1+\ldots+p_{s-1}x_{s-1}\equiv k\gcd(n,p_1,\ldots,p_{s-1}) (\modd\; n).
\end{equation}
Combining~\eqref{e:psk} and~\eqref{e:p1ps-1} we get the claim.
\kd

Let us now introduce the following definition that appeared in~\cite{SSB} (see also~\cite{But:10}).

\begin{definition}[Visualized Matrices]
{\rm A nonzero $A\in\Rpnn$ is called 
\begin{itemize}
\item[{\rm (i)}] {\em visualized} if $A_{i,j}\leq\lambda(A)$ for all $i,j$, and
\item[{\rm (ii)}] {\em strictly visualized} if it is visualized and 
$A_{i,j}=\lambda(A)$ if and only if $(i,j)\in\crit(A)$.
\end{itemize}
}
\end{definition}

By~\eqref{e:lambdaAcirc} we have that $\lambda(A)=\max(a_0,a_1,\ldots,a_{n-1})$ for
$A=\Circulant(a_0,\ldots,a_{n-1})$, implying that $\lambda(A)=\max\limits_{i,j=1}^n A_{i,j}$ 
for any circulant $A$. That is, any circulant matrix is visualized.
We will now argue that it is also strictly visualized.

\begin{definition}[Threshold Digraphs]
\label{def:thresh}
{\rm Let $A\in\Rpnn$ and $h\in\Rp$. Define the {\em threshold digraph} of $A$ with respect to $h$
as the subgraph of $\digr(A)$ containing all edges $(i,j)$ with $A_{i,j}\geq h$, and all
nodes that are beginning and end nodes of those edges. Denote this threshold graph by $\digr(A,h)$.} 
\end{definition}

\if{
\begin{lemma}
\label{l:circcirc}
If $A$ is a circulant matrix and $(i,j)\in\digr(A,h)$ then $\digr(A,h)$ contains
a cycle with edge $(i,j)$.
\end{lemma}
{\it Proof:} By Lemma~\ref{l:elementary}, $(i,j)$ belongs to a cycle $(i_1,\ldots,i_n)$ with 
with all edge weights of that cycle equal to the weight of $(i,j)$,
hence the claim follows.
\kd
}\fi

\begin{proposition}
\label{p:strictvis}
Let $A\in\Rpnn$ be a nonzero circulant matrix. Then it is strictly visualised,
and $\crit(A)=\digr(A,\lambda(A))$.
\end{proposition}
{\it Proof:}
By~\eqref{e:lambdaAcirc} no entry of $A$ exceeds $\lambda(A)$. Hence $A$ is 
visualized. Also recall that $\lambda(A)>0$ by Corollary~\ref{l:zerocase}.

If $A_{i,j}<\lambda(A)$ then the mean weight of any cycle with edge $(i,j)$ is strictly less than $\lambda(A)$,
so $(i,j)$ is not critical. In other words, $(i,j)$ being critical implies $A_{i,j}=\lambda(A)$.

It remains to show that if $A_{i,j}=\lambda(A)$, which is equivalent to $(i,j)$ being an edge of 
$\digr(A,\lambda(A))$, then $(i,j)$ is critical. In this case 
by Lemma~\ref{l:elementary} $(i,j)$ lies in a cycle with all edge weights equal to $\lambda(A)$.
The weights of all edges in this cycle are equal to $\lambda(A)$, hence the 
mean weight of this cycle is $\lambda(A)$, i.e., it is a critical cycle and $(i,j)$ is critical.
This completes the proof.
\kd

{\it Proof of Proposition~\ref{t:critcirc}.} 
First observe that Proposition~\ref{p:strictvis} implies that $\crit(A)=\digr(A,\lambda(A))$
and hence the set of critical edges of a circulant matrix $A$ is given by
\begin{equation}
\label{e:critedges}
E_c(A)=\{(i,j)\colon i=j\; \text{if $a_0=\lambda(A)$ or $j-i\equiv p_k(\modd\; n)$, $k\in\{1,\ldots,s\}$}\}
\end{equation}
where $p_1,\ldots,p_s$ are such that $a_{p_1}=\ldots=a_{p_s}=\lambda(A)$
(and $p_1>p_2>\ldots>p_s$).

We now consider the component of $\crit(A)$ which contains node $i$, for 
$i$ from the set $\{1,\ldots,\gcd(n,p_1,\ldots,p_s)\}$.

Let us argue that the node set of this component is given by
\begin{equation}
\label{e:onecomp}
\{k\in N\colon k\equiv i+l_1p_1+\ldots+l_sp_s(\modd\; n),\ l_1,\ldots,l_s\in\Nat\cup\{0\}\},
\end{equation}
Indeed, by~\eqref{e:critedges} edges 
$(i,j)$ where $j\equiv (l+p_t)(\modd\; n))$ for some $t\in\{1,\ldots, s\}$ are the only edges that issue from $i$
and are critical. Using this observation, the claim follows by simple induction.

Using Lemma~\ref{l:NT} we now observe that~\eqref{e:onecomp} is the same as
\begin{equation*}
\{i+k\gcd(n,p_1,\ldots,p_s))\colon k\in\{0,\ldots, (n/\gcd(n,p_1,\ldots,p_s))-1\}. 
\end{equation*}
This set does
not intersect with the node set of any component containing a different node in $\{1,\ldots,\gcd(n,p_1,\ldots,p_s)\}$, and this yields $\gcd(n,p_1,\ldots,p_s)$ strongly connected components of $\crit(A)$.
Isomorphism between two components containing $i_1\in\{1,\ldots,\gcd(n,p_1,\ldots,p_s)\}$
and $i_2\in\{1,\ldots,\gcd(n,p_1,\ldots,p_s)\}$ is induced by the following mapping on their set of nodes:
$$i_1+k\gcd(n,p_1,\ldots,p_s))\mapsto i_2+k\gcd(n,p_1,\ldots,p_s)).$$
This completes the proof of part (i) of Proposition~\ref{t:critcirc}.

If $a_0=\lambda(A)$ then $\crit(A)$ contains all loops of the form $(i,i)$ for $1\leq i\leq n$,
and the cyclicity of every component of $\crit(A)$ is $1$ since it contains a loop. 
When $a_0<\lambda(A)$, we can use the result of \cite{Gav-99}~Theorem~3.3 part (i) since 
this result describes the cyclicity of any component of the threshold
digraph $\digr(A,\lambda(A))$ (see \cite{Gav-99}~Theorem~3.1.),
and since $\crit(A)=\digr(A,\lambda(A))$ by Proposition~\ref{p:strictvis}.
According to this result, that cyclicity is equal to
any of the three expressions given in~\eqref{e:circcycl}.
This completes the proof of part (ii).
\kd

\subsection{Inclusion of attraction cones: Proof of Theorem~\ref{t:attrincl}}

\label{ss:proof2}

Before considering the problem of our interest, let us recall the notion of cyclic classes
which will be necessary for some proofs.
\begin{definition}[Cyclic Classes]
\label{def:cyclasses}
{\rm Let $\digr=(N,E)$ be a strongly connected graph with cyclicity 
$\sigma(\digr)$, and let $i,j\in N$. Nodes $i,j$ are said to belong to the 
same {\em cyclic class} if the lengths of some (and hence all) walks connecting $i$ to $j$
are a multiple of $\sigma(\digr)$. 

The cyclic class of $i$ will be denoted by $[i]$. We also write $[i]\to_1[j]$ if 
the lengths of some (and hence all) walks connecting a member of $[i]$ to a member of $[j]$
have length congruent to $1$ modulo $\sigma(\digr)$.

By cyclic classes of a completely reducible digraph we mean cyclic classes of
its (strongly connected) components.
}
\end{definition}

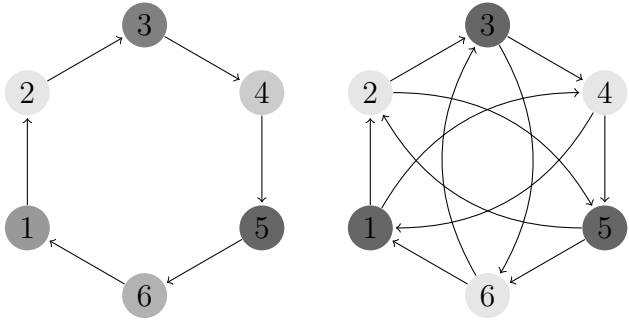
\begin{figure}
\begin{tabular}{ccc}
\begin{tikzpicture}[shorten >=1pt,->, scale=1.5]
  \tikzstyle{vertex1}=[circle,fill=black!10,minimum size=17pt,inner sep=1pt]
\tikzstyle{vertex2}=[circle,fill=black!50,minimum size=17pt,inner sep=1pt]
\tikzstyle{vertex3}=[circle,fill=black!20,minimum size=17pt,inner sep=1pt]
\tikzstyle{vertex4}=[circle,fill=black!60,minimum size=17pt,inner sep=1pt]
\tikzstyle{vertex5}=[circle,fill=black!30,minimum size=17pt,inner sep=1pt]
\tikzstyle{vertex6}=[circle,fill=black!40,minimum size=17pt,inner sep=1pt]


    \node[vertex1,xshift=0cm,yshift=0cm] (1) at (150:1.2cm) {$2$};
	\node[vertex2,xshift=0cm,yshift=0cm] (2) at (90:1.2cm) {$3$};
	\node[vertex3,xshift=0cm,yshift=0cm] (3) at (30:1.2cm) {$4$};
	\node[vertex4,xshift=0cm,yshift=0cm] (4) at (-30:1.2cm) {$5$};
	\node[vertex5,xshift=0cm,yshift=0cm] (5) at (-90:1.2cm) {$6$};
	\node[vertex6,xshift=0cm,yshift=0cm] (6) at (-150:1.2cm) {$1$};

\foreach \from/\to in {1/2,2/3,3/4,4/5,5/6,6/1}
    \draw (\from) -- (\to);
    
\end{tikzpicture} &&
\begin{tikzpicture}[shorten >=1pt,->, scale=1.5]
  \tikzstyle{vertex1}=[circle,fill=black!10,minimum size=17pt,inner sep=1pt]
\tikzstyle{vertex2}=[circle,fill=black!60,minimum size=17pt,inner sep=1pt]
\foreach \name/\angle/\text in {1/150/2, 
                                  3/30/4,  5/-90/6}
    \node[vertex1,xshift=0cm,yshift=0cm] (\name) at (\angle:1.2cm) {$\text$};
\foreach \name/\angle/\text in { 2/90/3, 
                                   4/-30/5,  6/-150/1}
    \node[vertex2,xshift=0cm,yshift=0cm] (\name) at (\angle:1.2cm) {$\text$};

\draw (6) to [bend left] (3);
\draw (3) to [bend left] (6);
\draw (1) to [bend left] (4);
\draw (4) to [bend left] (1);
\draw (2) to [bend left] (5);
\draw (5) to [bend left] (2);

\foreach \from/\to in {1/2,2/3,3/4,4/5,5/6,6/1}
    \draw (\from) -- (\to);
\end{tikzpicture}
\end{tabular}    
\caption{Cyclic classes of two graphs of Example~\ref{ex:twographs} (shown in different shades).
\label{fig:cyclic}}
\end{figure}

\begin{example}
\label{ex:twographs}
{\rm Consider two associated  graphs of $0$-$1$ matrices of Example~\ref{ex:01} shown in Figure~\ref{fig:cyclic}.
On the left, the graph consists just of one cycle of length $6$, hence its cyclicity is $6$
and the cyclic classes are $\{1\}$, $\{2\}$, $\{3\}$, $\{4\}$, $\{5\}$ and $\{6\}$. On the right,
the cyclicity of the graph is $2$ and the cyclic classes are $\{1,3,5\}$ and $\{2,4,6\}$.}\kd 
\end{example}

Cyclic classes are also called {\em components of imprimitivity}. We refer the reader to~\cite{BR} 
 Lemma 3.4.1
for a proof that belonging to the same cyclic class is a well-defined equivalence relation.

\begin{lemma}
\label{l:cyclasses}
Let $\digr$ be a strongly connected digraph.
\begin{itemize}
\item[{\rm (i)}] Let $\sigma(\digr)>1$ and let $i_0,i_1,\ldots,i_k$ be a walk on $\digr$.
Then $[i_{l-1}]\to_1[i_l]$ for each $l\in\{1,\ldots,k\}$.
\item[{\rm (ii)}] Let $C$ be a cycle of $\digr$. Then $C$ contains a member of 
each cyclic class of $\digr$. 
\end{itemize}
\end{lemma}
{\it Proof.} (i): Each edge is a walk of length $1$. Therefore $[i_{l-1}]\to_1[i_l]$ for each $l\in\{1,\ldots,k\}$.\\
(ii): Let $i$ be a node which is not in $C$. Let us show that $C$ contains a node in the cyclic class of $i$.
Since $\digr$ is strongly connected, there exists a walk connecting $i$ to a node $j$ of $C$. If the length of
this walk is a multiple of $\sigma(\digr)$ then $j\in [i]$. Otherwise, we concatenate this walk 
with a walk from $j$ to some node $k\in C$ whose edges belong to $C$ and such that the length of resulting walk is a multiple of $\sigma(\digr)$. Then $k\in [i]$ and the claim is proved.
\kd

We now derive a convenient form of a system defining the attraction cone 
for circulant matrices, based on the results of \cite{Ser-11}. 
Here $A^t_{i\bullet}$ denotes the $i$th row of $A^t$. We also write $i\sim_A j$ when
$i$ and $j$ belong to the same component of $\crit(A)$.

\begin{proposition}
\label{p:attrsys}
Let $A\in\Rpnn$ be a nonzero circulant matrix. Then
\begin{equation}
\label{e:attrsys1}
x\in\Attr(A)\Leftrightarrow  A^{n^2}_{i\bullet}\otimes x=A^{n^2}_{j\bullet}\otimes x\;\forall i,j\in N\,\text{s.t.}\, [i]\to_1 [j]
\end{equation}
and
\begin{equation}
\label{e:attrsys2}
x\in\Attr(A)\Leftrightarrow  A^{n^2}_{i\bullet}\otimes x=A^{n^2}_{j\bullet}\otimes x\; \forall i,j\in N\,\text{s.t.}\, i\sim_A j.
\end{equation}
\end{proposition}
{\it Proof.} 
By Proposition~\ref{p:attrdesc-circ}
\begin{equation}
\label{e:attrA1}
\Attr(A)=\{x\colon A^{n^2}\otimes x=A^{n^2+1}\otimes x\}
\end{equation}
\if{
By~\cite{Ser-11}~ Theorem 4.3, the system in~\eqref{e:attrA1} is 
equivalent to its critical subsystem and hence
\begin{equation*}
\Attr(A)=\{x\neq 0\colon A^t_{i\bullet}\otimes x=A^{t+1}_{i\bullet}\otimes x\quad\forall i\in N_c(A)\}
\end{equation*}
}\fi
Since $A$ is a circulant matrix, by Proposition~\ref{p:strictvis}
it is visualized, and then 
by~\cite{Ser-11} Proposition~2.8 we also have $A_{j\bullet}^{n^2}=A_{i\bullet}^{n^2+1}$ 
for any $i,j\in N_c(A)$ such that 
$[i]\to_1 [j]$.  This shows~\eqref{e:attrsys1}. To show~\eqref{e:attrsys2} recall that if a component of $\crit(A)$ has more than one cyclic class then for every two nodes $i,j$ of the component there is a walk $i_0=i,i_1,i_2,\ldots,i_k=j$ on $\crit(A)$ where $[i_{l-1}]\to_1 [i_l]$ for each $l\in\{1,\ldots,k\}$
by Lemma~\ref{l:cyclasses} part (i).
Hence $A^{n^2}_{i\bullet}\otimes x=A^{n^2}_{j\bullet}\otimes x$ holds for all nodes $i,j$ in that component.
If a component has only one cyclic class then~\cite{Ser-11} Proposition~2.8 implies that all rows with indices 
in that component are equal to each other, so the equations $A^{n^2}_{i\bullet}\otimes x=A^{n^2}_{j\bullet}\otimes x$
hold trivially for all pairs of nodes from that component.
\kd

It can be seen that we wrote out system~\eqref{e:attrsys1} for 
all examples of Section~\ref{s:circulants}. In the case of Example~\ref{ex:01}, 
for which $\digr(A)=\crit(A)$ and the cyclic classes are shown on Figure~\ref{fig:cyclic},
system~\eqref{e:attrsys1} reduces to $x_1=x_2=x_3=x_4=x_5=x_6$ for $A$ and to
$x_1\oplus x_3\oplus x_5=x_2\oplus x_4\oplus x_6$ for $B$.

We will also need the following observations.

\begin{lemma}
	\label{l:critincl}
Let $A,B\in\Rpnn$ be two matrices such that $\lambda(A)=\lambda(B)\neq 0$ and $A\leq B$. Then 
$\crit(A)\subseteq \crit(B)$. 
\end{lemma}
{\it Proof.}
Since $A\leq B$ the mean weight of each cycle in $B$ is not less than the mean weight of the same cycle in $A$. If that cycle is critical in $A$ then
	its mean weight $\lambda(A)$ cannot increase in $B$ since $\lambda(A)=\lambda(B)$. 
	Hence it equals $\lambda(B)$, i.e., the cycle belongs to $\crit(B)$.  
\kd  

\begin{lemma}
	\label{l:attrsysB}
Let $A,B\in\Rpnn$ be two circulant matrices with
$\lambda(A)=\lambda(B)\neq 0$, $A\leq B$. Then 
\begin{equation}
\label{e:attrsysB}
x\in\attr(B)\Leftrightarrow B^{n^2}_{i\bullet}\otimes x=B^{n^2}_{j\bullet}\otimes x
\quad\forall i,j\in N,\text{s.t.}\, i\sim_A j. 
\end{equation}
\end{lemma}
{\it Proof.}
By~\eqref{e:attrsys2},
\begin{equation}
\label{e:attrsys2B}
x\in \Attr(B)\Leftrightarrow B^{n^2}_{i\bullet}\otimes x=B^{n^2}_{j\bullet}\otimes x\; \forall i,j\in N\,\text{s.t.}\, i\sim_B j.
\end{equation}
We also have $\crit(A)\subseteq\crit(B)$ by Lemma~\ref{l:critincl} and hence
each $x\in\Attr(B)$ satisfies the system in~\eqref{e:attrsysB}. 

Suppose now that 
$x$ satisfies the system in\eqref{e:attrsysB}. We will show that $x$ also satisfies
\begin{equation}
\label{e:attrsysB1}
B^{n^2}_{i\bullet}\otimes x=B^{n^2}_{j\bullet}\otimes x\;\forall i,j\in N \,\text{s.t.}\, [i]\to_1 [j]
\end{equation}
so that $x\in\attr(B)$ by Proposition~\ref{p:attrsys}.
Since $\crit(A)\subseteq\crit(B)$, each component $\alpha$ of $\crit(A)$ belongs to a component $\beta$ of $\crit(B)$, and each component of $\crit(B)$ contains a component of $\crit(A)$ because $N_c(A)=N_c(B)=N$.
Hence it amounts to show that if $x$ satisfies the subsystem of equations in~\eqref{e:attrsysB}
corresponding to a component $\alpha$ of $\crit(A)$ then it also satisfies 
the subsystem of equations in~\eqref{e:attrsysB1} corresponding to the component $\beta$ of $\crit(B)$
such that $\alpha\subseteq\beta$.
But by Lemma~\ref{l:cyclasses} part (ii) 
each cyclic class of $\beta$ has a member in any cycle of $\beta$ and hence in any cycle of $\alpha$
(because $\alpha\subseteq\beta$).  This shows that for each $i,j$ with $[i]\to_1[j]$
in $\crit(B)$ there exist $k\in [i]$ and $l\in[j]$ on a cycle of $\alpha$ and then 
$B^{n^2}_{k\bullet}\otimes x= B^{n^2}_{l\bullet}\otimes x$ holds by~\eqref{e:attrsysB}. However,
$B^{n^2}_{k\bullet}=B^{n^2}_{i\bullet}$ and $B^{n^2}_{l\bullet}=B^{n^2}_{j\bullet}$ by \cite{Ser-11} Proposition~2.8.
Hence the claim follows.
\kd

Let us now introduce Kleene stars, as they will also be useful in the proof of Theorem~\ref{t:attrincl}.

\begin{definition}[Kleene Stars]
\label{def:kls}
{\rm Let $A\in\Rpnn$ have $\lambda(A)\leq 1$. Then
\begin{equation*}
A^*=I\oplus A\oplus A^2\oplus\ldots\oplus A^{n-1}
\end{equation*}
is called the {\em Kleene star of $A$}}.
\end{definition}


\begin{proposition}[e.g., \cite{But:10}, Corollary 1.6.16]
\label{p:kls}
Let $A\in\Rpnn$. Then $A^*=A$ if and only if one of the following
equivalent conditions hold: 
\begin{itemize}
\item[{\rm (i)}] $A^2=A$ and $A_{i,i}=1$ for all $i\in N$;
\item[{\rm (ii)}] $A_{i,i}=1$ and $A_{i,j}A_{j,k}\leq A_{i,k}$ for all $i,j,k\in N$.
\end{itemize}
\end{proposition}

More specifically, we will make use of the following.

\begin{lemma}
\label{l:n2kls}
Let $A\neq 0$ be a circulant matrix. Then 
$(A/\lambda(A))^{n^2}$ is a Kleene star.
\end{lemma}
{\it Proof.}
Note that $\lambda(A)\neq 0$ by Corollary~\ref{l:zerocase}.
By Proposition~\ref{p:kls} it suffices to show that $(A/\lambda(A))^{n^2}$ is an idempotent matrix and that $((A/\lambda(A))^{n^2})_{i,i}=1$
for all $i$. For the idempotency, observe that by Proposition~\ref{t:critcirc} part (ii) 
$\per(A)$ divides $n^2$, and that $T(A)\leq n^2$ by Proposition~\ref{p:transient-circ}. 
Hence  $(A/\lambda(A))^{2n^2}=(A/\lambda(A))^{n^2}$. 

For the remaining part of the claim, assume $\lambda(A)=1$ and recall that for any $t\geq 1$ 
and any $i,j\in N$, entry $(A^t)_{i,j}$ is equal to the greatest weight of 
a walk of length $t$ connecting $i$ to $j$ (e.g.,\cite{But:10}, Example 1.2.3).
Take $i\in\{1,\ldots,n\}$ and observe that
$\digr(A)$ contains a critical cycle of length $n$ going through $i$.
The weights of all entries of that cycle equal to $1$. Taking $n$ copies of this 
cycle we obtain a cycle in $\digr(A)$ of weight $1$ and length $n^2$.  The claim $((A/\lambda(A))^{n^2})_{i,i}=1$ follows since the weights of all entries and (therefore) of all walks are bounded by $1$.
\kd

We are now ready to prove the main result of Section~\ref{s:circulants}.

The proof will also make use of the following notation.
\begin{definition}
\label{def:moddodd}
{\rm Denote by $k[\modd\; n]$, respectively by $k[\modd'\; n]$,
the only number in $N_0=\{0,\ldots,n-1\}$, respectively in
$N=\{1,\ldots, n\}$, which is congruent to $k$ modulo $n$.}
\end{definition}

{\it Proof of Theorem~\ref{t:attrincl}.}
The case $\lambda(A)=\lambda(B)=0$ is trivial since in that case 
$A=B=0$ by Corollary~\ref{l:zerocase} and hence $\attr(A)=\attr(B)=\Rpn$. Otherwise,
as $\attr(A/\lambda(A))=\attr(A)$ and $\attr(B/\lambda(B))=\attr(B)$ (which follows, e.g., 
from Proposition~\ref{p:attrdesc}), we can 
assume without loss of generality $\lambda(A)=\lambda(B)=1$ and consider matrices $C=A^{n^2}$ and $D=B^{n^2}$. By Proposition~\ref{p:classical} 
$C$ and $D$ are circulants, hence 
$C=\Circulant(c_0,\ldots,c_{n-1})$ and 
$D=\Circulant(d_0,\ldots,d_{n-1})$ for some $c_0,\ldots,c_{n-1}$ and
$d_0,\ldots,d_{n-1}.$  Using $A\oplus B=B$ and the expansion for 
$(A\oplus B)^{n^2}$ we obtain $A^{n^2}\leq (A\oplus B)^{n^2}=B^{n^2}$, thus $C\leq D$. By Lemma~\ref{l:n2kls} 
   both of them are also Kleene stars. By Proposition~\ref{p:kls} we have 
	$D_{1,(\alpha+\gamma)[\modd'\; n]}\geq 
	D_{1,\alpha}\cdot D_{\alpha,(\alpha+\gamma)[\modd'\; n]}$ and hence
   	\begin{equation}
   	\label{e:dkls}
   	d_{(\alpha+\gamma-1)[\modd\; n]}\geq d_{\alpha-1}\cdot d_{\gamma}
   	\end{equation}
for any $\alpha\in\{1,\ldots,n\}$ and $\gamma\in\{0,\ldots,n-1\}$. In what follows we are going to show that the assumption that $\Attr(A)\subseteq \Attr(B)$ does not hold leads to a contradiction with~\eqref{e:dkls}
for some $\alpha$ and $\gamma$. 

 By Lemma~\ref{l:critincl} we have $\crit(A)\subseteq\crit(B)$.
 By Proposition~\ref{t:critcirc} part (i), 
 $\crit(A)$ consists of $l$ components whose node sets are of the form
 \begin{equation}
 \label{e:circcycles}
 \{k,k+l,k+2l,\ldots,k+(n/l-1)l\} \quad\text{for $k\in\{1,\ldots,l\}$},
 \end{equation}
 where $l$ is a divisor of $n$. Each of these node sets
belongs to some component of $\crit(B)$.


By Proposition~\ref{p:attrsys} $x\in\Attr(A)$ if and only if
\begin{equation}
\label{e:Csystem}
C_{k\bullet}\otimes x=C_{k+l\bullet}\otimes x=\ldots=C_{k+(n/l-1)l\bullet}\otimes x \quad\text{for $k\in\{1,\ldots,l\}$},
\end{equation}
and by Lemma~\ref{l:attrsysB} $x\in\Attr(B)$ if and only if
\begin{equation}
\label{e:Dsystem}
D_{k\bullet}\otimes x=D_{k+l\bullet}\otimes x=\ldots=D_{k+(n/l-1)l\bullet}\otimes x \quad\text{for $k\in\{1,\ldots,l\}$},
\end{equation}
We will refer to~\eqref{e:Csystem} or~\eqref{e:Dsystem} for fixed $k$ as to a {\em chain of equations}.

Suppose by contradiction that $x\in\Attr(A)$ but $x\notin\Attr(B)$. The latter means that there exist 
$k$ and $s$ such that
$D_{k\bullet}\otimes x> D_{k+ls\bullet}\otimes x$ for some integers $k$ and $s$. Assume without loss of generality that $k=1$ then
\begin{equation*}
D_{1\bullet}\otimes x=d_0 x_1\oplus d_1 x_2\oplus\ldots\oplus d_{n-1} x_n.
\end{equation*}

Let $d_{\alpha-1}\cdot x_{\alpha}$ be one of the terms where the maximum in the above expression is attained. In 
$D_{1+ls\bullet}\otimes x$ we find a term $d_{\alpha-1}\cdot x_{\beta}$ where $\alpha\equiv\beta(\modd\; l)$, and we have the inequality
$d_{\alpha-1}\cdot x_{\alpha}>d_{\alpha-1}\cdot x_{\beta}$ and hence $x_{\alpha}>x_{\beta}$.

Observe that $c_0=d_0=1$ since $C$ and $D$ are Kleene stars.
Since $\alpha\equiv\beta(\modd\; l)$ there exists a chain of equations among those of~\eqref{e:Csystem},
which contains both $c_0 x_{\alpha}=x_{\alpha}$ and $c_0 x_{\beta}=x_{\beta}$. The corresponding chain of equations holds (since
$x\in\attr(A)$), but $x_{\alpha}>x_{\beta}$ and therefore in the expression containing $c_0 x_{\beta}$ there is a term $c_{\gamma} x_{(\beta+\gamma)[\modd'\; n]}$ (for some $\gamma$) such that
$c_{\gamma} x_{(\beta+\gamma)[\modd'\; n]}\geq x_\alpha>0$, and hence
\begin{equation}
\label{e:firstineq}
d_{\gamma} x_{(\beta+\gamma)[\modd'\; n]}\geq x_\alpha. 
\end{equation}

Going back to the terms in the inequality $D_{1,\bullet}x> D_{1+ls,\bullet}x$ and knowing that the maximum in $D_{1,\bullet}x$ is attained at $d_{\alpha-1}x_{\alpha}$ and 
$D_{1+ls,\bullet}x$ contains a term of the form $d_{\alpha-1}x_{\beta}$, we see that $D_{1+ls,\bullet}x$ also contains the term 
$d_{(\alpha+\gamma-1)[\modd\; n]} x_{(\beta+\gamma)[\modd'\; n]}$ and that	
\begin{equation}
\label{e:secondineq}
d_{\alpha-1} x_{\alpha}>d_{(\alpha+\gamma-1)[\modd\; n]} x_{(\beta+\gamma)[\modd'\; n]}.
\end{equation}

Multiplying~\eqref{e:firstineq} by $d_{\alpha-1}$, combining with~\eqref{e:secondineq}
and canceling $x_{(\beta+\gamma)[\modd'\; n]}>0$ we have 
$$
d_{\alpha-1} d_{\gamma}>d_{(\alpha+\gamma-1)[\modd\; n]},
$$
which contradicts with the Kleene star property~\eqref{e:dkls}. The proof is complete.
\kd

\section{Acknowledgement}
We are grateful to the referees for their careful reading as well as 
numerous questions and helpful comments on the initial version of the paper, and we also would like to 
thank Dr. Michelle Delcourt for some useful suggestions regarding Abstract and Introduction.


\end{document}